\documentclass[12pt]{article}
\usepackage[english]{babel}
\usepackage{amssymb}
\usepackage{amsmath}
\usepackage{graphicx}
\usepackage{bbm}
\usepackage[right]{eurosym}
\usepackage[labelfont={bf,sf},font={small},
labelsep=space]{caption}
\usepackage{chngcntr}
\usepackage{xcolor}

\def\vs{\vspace}

\def\IN{\mathbb N}

\def\IR{\mathbb R}

\def\IQ{\mathbb Q}

\def\exp{\mathrm{exp}}

\def\ma{\mathcal}

\usepackage{filecontents}

\begin{document}

\begin{center}
{\Large\bf Restricted Log-Exp-Analytic Functions and some Differentiability Results}

\vs{0.3cm}
Andre Opris
\end{center}

\vs{0.2cm}
{\small {\bf Abstract.} In this article we define restricted log-exp-analytic functions as compositions of log-analytic functions and exponentials of locally bounded functions. We prove that the derivative of a restricted log-exp-analytic function is again restricted log-exp-analytic and that such a function exhibits strong quasianalytic properties. We establish the parametric version of Tamm’s theorem for this class of functions.}

\vs{0.7cm}
{\Large\bf Introduction}

\vs{0.2cm}
This paper contributes to analysis in the framework of o-minimal structures. O-minimality is a concept from mathematical logic with connections and applications to geometry, analysis, number theory and other areas. Sets and functions definable in an o-minimal structure (i.e. ‘belonging to’) exhibit tame geometric and combinatorial behaviour. We refer to the book of Van den Dries \cite{3} for the general properties of o-minimal structures (compare for example with \cite{3}, Chapter 3 for the notion of a definable cell decomposition); in the preliminary section we state the definition and give examples.

\vs{0.1cm}
By a deep model theoretical result of Van den Dries, Macintyre and Marker from \cite{4} every $\mathbb{R}_{\textnormal{an,exp}}$-definable function can be written as compositions of log-analytic functions and the global exponential ($\mathbb{R}_{\textnormal{an,exp}}$ is the o-minimal structure generated by all restricted analytic functions and the global exponential function, see \cite{2})). Log-analytic functions are iterated compositions from either side of globally subanalytic functions and the global logarithm (see \cite{7} respectively \cite{10} for elementary properties). They do not exhibit properties of the function $e^{-1/x}$ as flatness or infinite differentiability but not real analyticity (see \cite{7}). In this paper we define restricted log-exp-analytic functions as compositions from either side of log-analytic functions and exponentials of locally bounded functions and generalize the results from \cite{7} to this class of functions. 

\vs{0.4cm}
{\bf Example}

\vs{0.1cm}
The function
$$g: \textnormal{}]0,1[^2 \to \mathbb{R}, (t,x) \mapsto \arctan(\log(e^{1/t \cdot \log^2(1/x)}+\log(e^{e^{1/t}}+2))),$$ 
is restricted log-exp-analytic.

\vs{0.5cm}
Since the global logarithm is not globally subanalytic the class of restricted log-exp-analytic functions contains properly the class of $\mathbb{R}_{\textnormal{an}}$-definable functions. Although the global logarithm and the global exponential function are restricted log-exp-analytic the class of restricted log-exp-analytic functions is not a class of definable functions. It is properly contained in the class of $\mathbb{R}_{\textnormal{an,exp}}$-definable functions.

\vs{0.4cm}
{\bf Example}

\vs{0.1cm}
The function
$$h: \mathbb{R} \to \mathbb{R}, x \mapsto \left\{\begin{array}{ll} e^{-1/x}, & x > 0, \\
0, & \textnormal{ else, } \end{array}\right.$$
is not restricted log-exp-analytic, since one has $h \sim e^{-1/x}$ as $x \searrow 0$, but the function $h|_{\mathbb{R}_{>0}}$ is.

\vs{0.2cm}

\vs{0.1cm}
Since we give parametric versions for our main results we set up the concept of restricted log-exp-analytic functions in the variables $(w_1,...,w_l,u_1,...,u_m,x)$, where $(u_1,...,u_m,x)$ is serving as the tuple of independent variables of families of functions parameterized by $w:=(w_1,...,w_l)$. (The variable $x$ is needed to describe a preparation theorem for restricted log-exp-analytic functions with respect to a single variable which is suitable for our purposes.)

\vs{0.2cm}
Let $X \subset \mathbb{R}^l \times \mathbb{R}^m \times \mathbb{R}$ be $\mathbb{R}_{\textnormal{an,exp}}$-definable such that $X_w:=\{(u,x) \in \mathbb{R}^m \times \mathbb{R} \mid (w,u,x) \in X\}$ is open for every $w \in \mathbb{R}^l$. Call a function $f:X \to \mathbb{R}, (w,u,x) \mapsto f(w,u,x),$ restricted log-exp-analytic in $(u,x)$ if $f$ is the composition from either side of log-analytic functions and exponentials of locally bounded functions in $(u,x)$ ($g:X \to \mathbb{R}$ is locally bounded in $(u,x)$ if $g_w:X_w \to \mathbb{R}, (u,x) \mapsto g(w,u,x),$ is locally bounded for every $w \in \mathbb{R}^l$).

\vs{0.1cm}
As in the log-analytic case in \cite{7} a sufficiently good representation for such a restricted log-exp-analytic function can be obtained with preparation theorems for $\mathbb{R}_{\textnormal{an,exp}}$-definable functions (compare with Lion-Rolin \cite{8} for original versions, Van den Dries and Speissegger \cite{6} for a preparation in a more general context, Opris \cite{10} for strong versions and Cluckers and Miller \cite{1} for another kind of preparation for the special class of constructible functions suitable for questions on integration). 

\vs{0.1cm}
Outgoing from a preparation theorem for $\mathbb{R}_{\textnormal{an,exp}}$-definable functions from \cite{10} adapted to the restricted log-exp-analytic case our initial result is the following key observation for a restricted log-exp-analytic function $f:X \to \mathbb{R}, (w,u,x) \mapsto f(w,u,x),$ in $(u,x)$ where $X \subset \mathbb{R}^l \times \mathbb{R}^m \times \mathbb{R}$ is $\mathbb{R}_{\textnormal{an,exp}}$-definable such that $X_w$ is open and $0$ is interior point of $X_{(w,u)}$ for every $(w,u) \in \pi(X)$ (here $\pi:\mathbb{R}^l \times \mathbb{R}^m \times \mathbb{R} \to \mathbb{R}^l \times \mathbb{R}^m, (w,u,x) \mapsto (w,u),$ denotes the projection on the first $l+m$ coordinates).

\vs{0.1cm}
There is an $\mathbb{R}_{\textnormal{an,exp}}$-definable cell decomposition $\mathcal{C}$ of $X$ such that for every simple $C \in \mathcal{C}$ (see Definition 3.4) we have
$$f(t,x)=a(t)\vert{y_0(t,x)}\vert^{q_0} \cdot ... \cdot \vert{y_r(t,x)}\vert^{q_r}u(t,x)$$
for $(t,x) \in C$ where $y_0(t,x)=x$, $y_1(t,x)=\log(x)$, $y_2(t,x)=\log(-\log(x))$,
..., $y_r(t,x)=\log_{r-1}(-\log(x))$, the $q_j$'s are rational exponents, the function $a: \pi(C) \to \mathbb{R}, (w,u) \mapsto a(w,u),$ can be extended to a restricted log-exp-analytic function $\tilde{a}:\pi(X) \to \mathbb{R},(u,w) \mapsto \tilde{a}(u,w),$ in $u$ (i.e. $\tilde{a}|_{\pi(C)} = a$) and $u(t,x)$ is a unit of a special form which we describe below (where $t:=(w,u)$). 

\vs{0.1cm}
An immediate consequence is that the class of restricted log-exp-analytic functions is closed under taking derivatives. 

\vs{0.4cm}
{\bf Theorem A}

\vs{0.1cm}
{\it
Let $X \subset \mathbb{R}^n \times \mathbb{R}^m$ be $\mathbb{R}_{\textnormal{an,exp}}$-definable such that $X_t$ is open for every $t \in \mathbb{R}^n$. Let $f:X \to \mathbb{R}, (t,u) \mapsto f(t,u),$ be restricted log-exp-analytic in $u$. Let $i \in \{1,...,m\}$ be such that $f$ is differentiable with respect to $u_i$ on $X$. Then $\partial f/\partial u_i$ is restricted log-exp-analytic in $u$.} 

\vs{0.3cm}
One gets also a result on strong quasianalyticity (see Miller \cite{9} for this result in polynomially bounded o-minimal structures).

\vs{0.4cm}
{\bf Theorem B}

\vs{0.1cm}
{\it
Let $X \subset \mathbb{R}^n \times \mathbb{R}^m$ be definable such that $X_t$ is open and connected for every $t \in \mathbb{R}^n$. Let $f:X \to \mathbb{R}, (t,u) \mapsto f(t,u),$ be restricted log-exp-analytic in $u$. Then there is $N \in \mathbb{N}$ with the following property. If $f(t,-)$ is $C^N$ for $t \in \mathbb{R}^n$ and if there is $a \in X_t$ such that all derivatives up to order $N$ vanish in $a$ then $f(t,-)$ vanishes identically.} 

\vs{0.2cm}
By adapting the arguments of Van den Dries and Miller from \cite{5} to the restricted log-exp-analytic case one obtains the parametric version of Tamm's theorem for restricted log-exp-analytic functions.

\vs{0.4cm}
{\bf Theorem C}

\vs{0.1cm}
{\it
Let $X \subset \mathbb{R}^n \times \mathbb{R}^m$ be definable such that $X_t$ is open for $t \in \mathbb{R}^n$. Let $f:X \to \mathbb{R}, (t,u) \mapsto f(t,u),$ be restricted log-exp-analytic in $u$. Then there is $N \in \mathbb{N}$ such that the following holds for every $(t,u) \in X$: If $f(t,-)$ in $C^N$ at $u$ then $f(t,-)$ is real analytic at $u$.} 

\vs{0.3cm}
(Compare also with \cite{5} for a formulation of Theorem C in the globally subanalytic setting and with \cite{11} for Tamm's original version. Compare also with \cite{7} for a formulation of Theorem A, Theorem B and Theorem C in the log-analytic setting.)

\vs{0.2cm}
Theorem C implies that the property of being real analytic is a definable property in case of restricted log-exp-analytic functions. This shows that this class of functions shares its properties from the viewpoint of analysis with log-analytic and even with globally subanalytic functions. But this is not true for $\mathbb{R}_{\textnormal{an,exp}}$-definable functions in general as remarked in \cite{7} and in the end of \cite{5}.

\vs{0.5cm}
{\bf Example} (Kaiser/Opris, \cite{7})

\vs{0.1cm}
Consider the function 
$$f: \mathbb{R} \times \mathbb{R} \to \mathbb{R},
(t,x)\mapsto 
\left\{\begin{array}{lll}| {x} |^{|2t|},&& x \neq 0, \\
&\mbox{if}&\\
0,&& x=0,
\end{array}\right.$$
which is definable in $\mathbb{R}_{\exp}$. Then the set of all $t \in \mathbb{R}$ such that $f(t,-)$ is real-analytic at $0$ is the set of integers. 

\vs{0.5cm}
The paper is organized as follows. After a preliminary section on o-minimality, notations and conventions we give in Section 1 the definition of log-analytic functions, the exponential number and give a preparation theorem for definable functions. In Section 2 we introduce restricted log-exp-analytic functions, adapt the preparation theorem to this class of functions and give a version on simple sets. Section 3 is devoted to the proof of the above theorems.

\vs{0.8cm}
{\Large\bf Preliminaries}

\vs{0.3cm}
{\bf Semialgebraic sets:}

\vs{0.1cm}
A subset $A$ of $\mathbb{R}^n$, $n \geq 1$, is called \textbf{semialgebraic} if there are $k,l \in \mathbb{N}_0$ and real polynomials $f_i,g_{i,1},...,g_{i,k} \in \mathbb{R}[X_1,...,X_n]$ for $1 \leq i \leq l$ such that 
$$A = \bigcup_{i=1}^l \{x \in \mathbb{R}^n \mid f_i(x)=0, g_{i,1}(x)>0,...,g_{i,k}(x)>0\}.$$
A map is called semialgebraic if its graph is semialgebraic.

\vs{0.3cm}
{\bf Semi- and subanalytic sets:}

\vs{0.1cm}
A subset $A$ of $\mathbb{R}^n$, $n \geq 1$, is called \textbf{semianalytic} if for each $a \in \mathbb{R}^n$ there are open neighbourhoods $U,V$ of $a$ with $\overline{U} \subset V$, $k,l \in \mathbb{N}_0$ and real analytic functions $f_i,g_{i,1},...,g_{i,k}$ on $V$ for $1 \leq i \leq l$, such that
$$A \cap U = \bigcup_{i=1}^l\{x \in U \mid f_i(x)=0, g_{i,1}(x)>0,...,g_{i,k}(x)>0\}.$$
A subset $B$ of $\mathbb{R}^n$, $n \geq 1$, is called \textbf{subanalytic} if for each $a \in \mathbb{R}^n$ there is an open neighbourhood $U$ of $a$, some $p \geq n$ and some bounded semianalytic set $A \subset \mathbb{R}^p$ such that $B \cap U=\pi_n(A)$ where $\pi_n:\mathbb{R}^p \to \mathbb{R}^n, (x_1,...,x_p) \mapsto (x_1,...,x_n),$ is the projection on the first $n$ coordinates.
A map is called semianalytic or subanalytic if its graph is a semianalytic resp. subanalytic set, respectively. A set is called \textbf{globally semianalytic} or \textbf{globally subanalytic} if it is semianalytic or subanalytic, respectively, in the ambient projective space (or equivalently, after applying the semialgebraic homeomorphism $\mathbb{R}^n \to \textnormal{}]-1,1[^n, x_i \mapsto x_i/\sqrt{1+x_i^2}$.)

\vs{0.3cm}
{\bf O-minimal structures:}

\vs{0.1cm}
A \textbf{structure} on $\mathbb{R}$ is axiomatically defined as follows. For $n \in \mathbb{N}$ let $M_n$ be a set of subsets of $\mathbb{R}^n$ and let $\mathcal{M}:=(M_n)_{n \in \mathbb{N}}$. Then $\mathcal{M}$ is a structure on $\mathbb{R}$ if the following holds for all $m,n,p \in \mathbb{N}$.
\begin{itemize}
	\item [(S1)] If $A,B \in M_n$ then $A \cup B$, $A \cap B$ and $\mathbb{R}^n \setminus A \in M_n$. (So $M_n$ is a Boolean algebra of subsets of $M_n$.)
	\item [(S2)] If $A \in M_n$ and $B \in M_m$ then $A \times B \in M_{n+m}$.
	\item [(S3)] If $A \in M_p$ and $p \geq n$ then $\pi_n(A) \in M_n$ where $\pi_n:\mathbb{R}^p \to \mathbb{R}^n, (x_1,...,x_p) \mapsto (x_1,...,x_n)$, denotes the projection on the first $n$ coordinates.
	\item [(S4)] $M_n$ contains the semialgebraic subsets of $\mathbb{R}^n$. 
\end{itemize}
The structure $\mathcal{M}=(M_n)_{n \in \mathbb{N}}$ on $\mathbb{R}$ is called \textbf{o-minimal} if additionally the following holds.
\begin{itemize}
	\item [(O)] The sets in $M_1$ are exactly the finite unions of intervals and points.
\end{itemize}

A subset of $\mathbb{R}^n$ is called \textbf{definable} in the structure $\mathcal{M}$ if it belongs to $M_n$.
A function is definable in $\mathcal{M}$ if its graph is definable in $\mathcal{M}$.
The o-minimality axiom (O) implies that a subset of $\mathbb{R}$ which is definable in an o-minimal structure on $\mathbb{R}$ has only finitely many connected components. But much more can be deduced from o-minimality. A definable subset of $\mathbb{R}^n$, $n \in \mathbb{N}$ arbitrary, has only finitely many connected components and these are again definable. More generally, sets and functions definable in an o-minimal
structure exhibit tame geometric behaviour, for example the existence of definable
cell decomposition. We refer to the book of Van den Dries \cite{3} for this and more of the general properties of o-minimal structures.

\vs{0.5cm}
{\bf Examples of o-minimal structures:}

\begin{itemize}
	\item [(1)] The smallest o-minimal structure on $\mathbb{R}$ is given by the semialgebraic sets. It is denoted by $\mathbb{R}$.
	\item [(2)] $\mathbb{R}_{\textnormal{exp}}$, the structure generated on the real field by the global exponential function $\exp:\mathbb{R} \to \mathbb{R}_{>0}$ (i.e. the smallest structure containing the semialgebraic sets and the graph of the exponential function), is o-minimal.
	\item [(3)] $\mathbb{R}_{\textnormal{an}}$, the structure generated on the real field by the restricted analytic functions, is o-minimal. A function $f:\mathbb{R}^n \to \mathbb{R}$ is called restricted analytic if there is a function $g$ that is real analytic on a neighbourhood of $[-1,1]^n$ such that $f=g$ on $[-1,1]^n$ and $f=0$ else. The sets definable in $\mathbb{R}_{\textnormal{an}}$ are precisely the globally subanalytic ones.
	\item [(4)] $\mathbb{R}_{\textnormal{an,exp}}$, the structure generated by $\mathbb{R}_{\textnormal{an}}$ and $\mathbb{R}_{\textnormal{exp}}$, is o-minimal.
\end{itemize}

\vs{0.2cm}
{\bf Notation}

\vs{0.1cm}
The empty sum is by definition $0$ and the empty product is by definition $1$. By $\mathbb{N}=\{1,2,...\}$ we denote the set of natural numbers and by $\mathbb{N}_0=\{0,1,2,...\}$ the set of nonnegative integers. Given $x \in \mathbb{R}$ let $\lceil{x}\rceil$ be the smallest integer which is not smaller than $x$. We set $\mathbb{R}_{>0}:=\{x \in \mathbb{R} \mid x > 0\}$. For $m,n \in \mathbb{N}$ we denote by $M(m \times n,\mathbb{Q})$ the set of $m \times n$-matrices with rational entries and for $P \in M(m \times n,\mathbb{Q})$ we denote by $^tP \in M(n \times m,\mathbb{Q})$ its transpose. Given $x \in \mathbb{R} \setminus \{0\}$ let $\textnormal{sign}(x) \in \{\pm 1\}$ its sign. We use the usual $o$-notation and $O$-notation. By the symbol ''$\sim$'' we denote asymptotic equivalence. By $\log_k$ we denote the $k$-times iterated of the natural logarithm.

\vs{0.1cm}
For $m \in \mathbb{N}$ and a set $X \subset \mathbb{R}^m$ we set the following: For a set $E$ of positive real valued functions on $X$ we set $\log(E):=\{\log(g) \mid g \in E\}$.

\vs{0.2cm}
For $X \subset \mathbb{R}^n \times \mathbb{R}^m$ and $t \in \mathbb{R}^n$ we set $X_t:=\{x \in \mathbb{R}^m \mid (t,x) \in X\}$ and for a function $f:X \to \mathbb{R}$ and $t \in \mathbb{R}^n$ we set $f_t:X_t \to \mathbb{R}, x \mapsto f(t,x)$. 

\vs{0.2cm}
For $X \subset \mathbb{R}^n \times \mathbb{R}$ we set $X_{\neq 0}:=\{(t,x) \in X \mid x \neq 0\}$.

\vs{0.4cm}
{\bf Convention}

\vs{0.1cm}
Definable means definable in $\mathbb{R}_{\textnormal{an,exp}}$ if not otherwise mentioned.
 
\section{Definable Functions}

\subsection{Log-Analytic Functions and the \\ Exponential Number}

Compare also with \cite{10} for Section 1.1. 

\vs{0.3cm}
{\bf1.1 Definition} 

\vs{0.1cm}
Let $f:X \to \mathbb{R}$ be a function.
\begin{itemize}
	\item [(a)] Let $r \in \mathbb{N}_0$. By induction on $r$ we define that $f$ is \textbf{log-analytic of order at most} $r$.
	
	\vs{0.3cm}
	\textbf{Base case}: The function $f$ is log-analytic of order at most $0$ if $f$ is piecewise the restriction of  globally subanalytic functions, i.e. there is a decomposition $\mathcal{C}$ of $X$ into finitely many definable sets such that for $C \in \mathcal{C}$ there is a globally subanalytic function $F:\mathbb{R}^m \to \mathbb{R}$ such that $f|_C = F|_C$.
	
	\vs{0.3cm}
	\textbf{Inductive step}: The function $f$ is log-analytic of order at most $r$ if the following holds: There is a decomposition $\mathcal{C}$ of $X$ into finitely many definable sets such that for $C \in \mathcal{C}$ there are $k,l \in \mathbb{N}_{0}$, a globally subanalytic function $F:\mathbb{R}^{k+l} \to \mathbb{R}$, and log-analytic functions $g_1,...,g_k:C \to \mathbb{R}, h_1,...,h_l:C \to \mathbb{R}_{>0}$ of order at most $r-1$ such that
	$$f|_C=F(g_1,...,g_k,\log(h_1),...,\log(h_l)).$$
	
	\item[(b)] Let $r \in \mathbb{N}_0$. We call $f$ \textbf{log-analytic of order} $r$ if $f$ is log-analytic of order at most $r$ but not of order at most $r-1$.
	
	\item[(c)] We call $f$ \textbf{log-analytic} if $f$ is log-analytic of order $r$ for some $r \in \mathbb{N}_0$.
\end{itemize}

\vs{0.3cm}
{\bf1.2 Remark}
\begin{itemize}
	\item [(1)] Let $r \in \mathbb{N}_0$. The set of log-analytic functions on $X$ of order at most $r$ is an $\mathbb{R}$-algebra with respect to pointwise addition and multiplication.
	\item [(2)] The set of log-analytic functions on $X$ is an $\mathbb{R}$-algebra with respect to pointwise addition and multiplication.
\end{itemize}

\vs{0.3cm}
{\bf1.3 Remark}

\vs{0.1cm}
A definable function $f:X \to \mathbb{R}$ is positive if and only if there is a definable function $g:X \to \mathbb{R}$ such that $f(x)=\exp(g(x))$ for every $x \in X$.

\vs{0.5cm}
{\bf1.4 Definition}

\vs{0.1cm}
Let $f:X \to \mathbb{R}$ be a function. Let $E$ be a set of positive definable functions on $X$.
\begin{itemize}
	\item [(a)] By induction on $e \in \mathbb{N}_0$ we define that $f$ has \textbf{exponential number at most $e$ with respect to $E$}.
	
	\vs{0.2cm}
	{\bf Base Case}: The function $f$ has exponential number at most $0$ with respect to $E$ if $f$ is log-analytic.
	
	\vs{0.2cm}
	{\bf Inductive Step}: The function $f$ has exponential number at most $e$ with respect to $E$ if the following holds: There are $k,l \in \mathbb{N}_0$, functions $g_1,...,g_k:X \to \mathbb{R}$ and $h_1,...,h_l:X \to \mathbb{R}$ with exponential number at most $e-1$ with respect to $E$ and a log-analytic function $F:\mathbb{R}^{k+l} \to \mathbb{R}$ such that
	$$f=F(g_1,...,g_k,\exp(h_1),...,\exp(h_l))$$
	and $\exp(h_1),...,\exp(h_l) \in E$.
	
	\item [(b)] Let $e \in \mathbb{N}_0$. We say that $f$ has \textbf{exponential number $e$ with respect to $E$} if $f$ has exponential number at most $e$ with respect to $E$ but not at most $e-1$ with respect to $E$.
	
	\item [(c)] We say that $f$ \textbf{can be constructed from $E$} if there is $e \in \mathbb{N}_0$ such that $f$ has exponential number $e$ with respect to $E$. 
\end{itemize}

\vs{0.5cm}
{\bf1.5 Remark}

\vs{0.1cm}
Let $E$ be a set of positive definable functions on $X$. Let $f$ be a function on $X$ which can be constructed from $E$. The following holds.
\begin{itemize}
	\item[(1)] $f$ has a unique exponential number $e \in \mathbb{N}_0$ with respect to $E$ and is definable.
	\item[(2)] $f$ has exponential number $0$ with respect to $E$ if and only if it is log-analytic. 
\end{itemize}

By a quantifier elimination result such as Corollary 4.7 in \cite {4} we obtain the following.

\vs{0.3cm}
{\bf1.6 Remark}

\vs{0.1cm}
Let $E$ be the set of all positive definable functions on $X$. Then every definable function $f:X \to \mathbb{R}$ can be constructed from $E$.

\vs{0.5cm}
{\bf1.7 Remark}

\vs{0.1cm}
Let $X_1,X_2 \subset \mathbb{R}^m$ be definable and disjoint. Let $X = X_1 \cup X_2$. For $j \in \{1,2\}$ let $E_j$ be a set of positive definable functions on $X_j$ and $f_j:X_j \to \mathbb{R}$ be a function. Let $e \in \mathbb{N}_0$ be such that $f_j$ has exponential number at most $e$ with respect to $E_j$ for $j \in \{1,2\}$. Let 
$$E:=\{g \mid g:X \to \mathbb{R} \textnormal{ is a function with } g|_{X_j} \in E_j \textnormal{ for } j \in \{1,2\}\}.$$
Then 
$$f:X \to \mathbb{R}, x \mapsto \left\{\begin{array}{ll} f_1(x) , & x \in X_1,  \\
f_2(x), & x \in X_2, \end{array}\right.$$
has exponential number at most $e$ with respect to $E$.

\vs{0.5cm}
{\bf1.8 Remark}

\vs{0.1cm}
Let $e \in \mathbb{N}_0$. Let $E$ be a set of positive definable functions on $X$.
\begin{itemize}
	\item[(1)] Let $f:X \to \mathbb{R}$ be a function with exponential number at most $e$ with respect to $E$. Then $\exp(f)$ has exponential number at most $e+1$ with respect to $E \cup \{\exp(f)\}$.
	\item[(2)] Let $s \in \mathbb{N}_0$. Let $f_1,...,f_s:X \to \mathbb{R}$ be functions with exponential number at most $e$ with respect to $E$ and let $F:\mathbb{R}^s \to \mathbb{R}$ be log-analytic. Then $F(f_1,...,f_s)$ has exponential number at most $e$ with respect to $E$. 
\end{itemize}

\vs{0.3cm}
{\bf Proof}

\vs{0.1cm}
Compare with \cite{10}, Section 1. \hfill$\blacksquare$

\subsection{A Preparation Theorem for $\mathbb{R}_{\textnormal{an,exp}}$-definable \\
	Functions}

Let $n \in \mathbb{N}$. Let $t$ range over $\mathbb{R}^n$ and $x$ over $\mathbb{R}$. We fix a definable set $C \subset \mathbb{R}^n \times \mathbb{R}$. let $\pi: \mathbb{R}^n \times \mathbb{R} \to \mathbb{R}, (t,x) \mapsto t,$ be the projection on the first $n$ coordinates. 

\vs{0.5cm}
{\bf1.9 Definition}

\vs{0.1cm}
Let $r \in \mathbb{N}_0$. A tuple $\mathcal{Y}:=(y_0,...,y_r)$ of functions on $C$ is called an \textbf{$r$-logarithmic scale} on $C$ with \textbf{center} $\Theta=(\Theta_0,...,\Theta_r)$ if the following holds:
\begin{itemize}
	\item[(a)] $y_j>0$ or $y_j<0$ for every $j \in \{0,...,r\}$.
	\item[(b)] $\Theta_j$ is a definable function on $\pi(C)$ for every $j \in \{0,...,r\}$.
	\item[(c)] We have $y_0(t,x)=x-\Theta_0(t)$ and for every $j \in \{1,...,r\}$ we have $y_j(t,x)=\log(\vert{y_{j-1}(t,x)}\vert) - \Theta_j(t)$ for all $(t,x) \in C$.
	\item[(d)] Either there is $\epsilon_0 \in \textnormal{}]0,1[$ such that $0<\vert{y_0(t,x)}\vert < \epsilon_0\vert{x}\vert$ for all $(t,x) \in C$ or $\Theta_0=0$, and for every $j \in \{1,...,r\}$ either there is $\epsilon_j \in \textnormal{}]0,1[$ such that $0<\vert{y_j(t,x)}\vert<\epsilon_j\vert{\log(\vert{y_{j-1}(t,x)}\vert)}\vert$ for all $(t,x) \in C$ or $\Theta_j=0$.
\end{itemize}

\vs{0.5cm}
{\bf1.10 Definition}

\vs{0.1cm}
We call $g:\pi(C) \to \mathbb{R}$ a \textbf{$C$-heir} if there is $r \in \mathbb{N}_0$, an $r$-logarithmic scale $\mathcal{Y}$ with center $(\Theta_0,...,\Theta_r)$ on $C$, and $l \in \{1,...,r\}$ such that $g=\exp(\Theta_l)$.

\vs{0.5cm}
{\bf1.11 Definition}

\vs{0.1cm}
We call $g:\pi(C) \to \mathbb{R}$ \textbf{$C$-nice} if there is a set $E$ of $C$-heirs such that $g$ can be constructed from $E$.

\vs{0.5cm}
Note that the class of log-analytic functions on $\pi(C)$ is a proper subclass of the class of $C$-nice functions. (Compare with \cite{10} for examples and further elementary properties of $C$-heirs and $C$-nice functions.)

\vs{0.5cm}
{\bf1.12 Definition}

\vs{0.1cm}
Let $r \in \mathbb{N}_0$, $\mathcal{Y}:=(y_0,...,y_r)$ be a logarithmic scale on $C$ and $q:=(q_0,...,q_r) \in \mathbb{Q}^{r+1}$. For $(t,x) \in C$ we write 
$$\vert{\mathcal{Y}(t,x)}\vert^{\otimes q}:=\prod_{j=0}^r \vert{y_j(t,x)}\vert^{q_j}.$$

\vs{0.2cm}
{\bf1.13 Definition}

\vs{0.1cm}
Let $r \in \mathbb{N}_0$. Let $f:C \to \mathbb{R}$ be a function. We say that $f$ is \textbf{$r$-log-analytically prepared in $x$ with center $\Theta$} if
$$f(t,x)=a(t) \vert{\mathcal{Y}(t,x)}\vert^{\otimes q}u(t,x)$$
for all $(t,x) \in C$ where $a$ is a definable function on $\pi(C)$ which vanishes identically or has no zero, $\mathcal{Y}=(y_0,...,y_r)$ is an $r$-logarithmic scale with center $\Theta$ on $C$, $q \in \mathbb{Q}^{r+1}$ and the following holds for $u$. There is $s \in \mathbb{N}$ such that $u=v \circ \phi$ where $v$ is a power series which converges on an open neighbourhood of $[-1,1]^s$ with $v([-1,1]^s) \subset \mathbb{R}_{>0}$ and $\phi:=(\phi_1,...,\phi_s):C \to [-1,1]^s$ is a function of the form 
$$\phi_j(t,x):=b_j(t)\vert{\mathcal{Y}(t,x)}\vert^{\otimes p_j}$$
for $j \in \{1,...,s\}$ and $(t,x) \in C$ where $b_j:\pi(C) \to \mathbb{R}$ is definable for $j \in \{1,...,s\}$ and $p_j:=(p_{j0},...,p_{jr}) \in \mathbb{Q}^{r+1}$. We call $a$ a \textbf{coefficient} and $b:=(b_1,...,b_s)$ a tuple of \textbf{base functions} for $f$. An \textbf{LA-preparing tuple} for $f$ is then
$$\mathcal{J}:=(r,\mathcal{Y},a,q,s,v,b,P)$$
where
$$P:=\left(\begin{array}{cccc}
p_{10}&\cdot&\cdot&p_{1r}\\
\cdot&& &\cdot\\
\cdot&& &\cdot\\
p_{s0}&\cdot&\cdot&p_{sr}\\
\end{array}\right)\in M\big(s\times (r+1),\mathbb{Q}).$$

\vs{0.5cm}
{\bf1.14 Fact} (Opris, \cite{10}, Theorem A)

\vs{0.1cm}
{\it
Let $m \in \mathbb{N}$, $r \in \mathbb{N}_0$. Let $X \subset \mathbb{R}^n \times \mathbb{R}$ be definable. Let $f_1,....,f_m:X \to \mathbb{R}$ be log-analytic functions of order at most $r$. Then there is a definable cell decomposition $\mathcal{C}$ of $X_{\neq 0}$ such that $f_1|_C,...,f_m|_C$ are $r$-log-analytically prepared in $x$ with $C$-nice coefficients, $C$-nice base functions and common $C$-nice center for $C \in \mathcal{C}$.}

\vs{0.5cm}
{\bf1.15 Definition}

\vs{0.1cm}
Let $r \in \mathbb{N}_0$ and $e \in \mathbb{N}_0 \cup \{-1\}$. Let $C \subset \mathbb{R}^n \times \mathbb{R}$ be definable and $f:C \to \mathbb{R}$ be a function. Let $E$ be a set of positive definable functions on $C$.

\vs{0.1cm}
By induction on $e$ we define that $f$ is \textbf{(purely) $(e,r)$-prepared} in $x$ with center $\Theta$ with respect to $E$. To this preparation we assign a \textbf{(pure) preparing tuple} for $f$.

\vs{0.3cm}
$e=-1$: We say that $f$ is (purely) $(-1,r)$-prepared in $x$ with center $\Theta$ with respect to $E$ if $f$ is the zero function. A (pure) preparing tuple for $f$ is $(0)$.

\vs{0.3cm}
$e-1 \to e$: We say that $f$ is (purely) $(e,r)$-prepared in $x$ with center $\Theta$ with respect to $E$ if 
$$f(t,x)=a(t)\vert{\mathcal{Y}(t,x)}\vert^{\otimes q}\exp(c(t,x)) u(t,x)$$
for every $(t,x) \in C$ where $a:\pi(C) \to \mathbb{R}$ is $C$-nice (log-analytic), $\mathcal{Y}:=(y_0,...,y_r)$ is an $r$-logarithmic scale with $C$-nice (log-analytic) center $\Theta$, $q=(q_0,...,q_r) \in \mathbb{Q}^{r+1}$, $\exp(c) \in E$ where $c$ is (purely) $(e-1,r)$-prepared in $x$ with center $\Theta$ with respect to $E$ and the following holds for $u$. There is $s \in \mathbb{N}$ such that $u=v \circ \phi$ where $\phi:=(\phi_1,...,\phi_s):C \to [-1,1]^s$ with 
$$\phi_j(t,x)=b_j(t)\vert{y_0(t,x)}\vert^{q_0} \cdot ... \cdot \vert{y_r(t,x)}\vert^{q_r}\exp(d_j(t,x))$$
for $j \in \{1,...,s\}$ and $(t,x) \in C$, $b_j:\pi(C) \to \mathbb{R}$ is $C$-nice (log-analytic), $\exp(d_j) \in E$, $d_j:C \to \mathbb{R}$ is (purely) $(e-1,r)$-prepared in $x$ with center $\Theta$ with respect to $E$ and $v$ is a power series which converges absolutely on an open neighbourhood of $[-1,1]^s$ with $v([-1,1]^s) \subset \mathbb{R}_{>0}$.

\vs{0.5cm}
{\bf1.16 Remark}

\vs{0.1cm}
Let $e \in \mathbb{N}_0 \cup \{-1\}$ and $r \in \mathbb{N}_0$. Let $E$ be a set of positive definable functions on $C$. 
\begin{itemize}
	\item[(1)]
	Let $f:C \to \mathbb{R}$ be $(e,r)$-prepared in $x$ with respect to $E$. Then there is a set $\mathcal{E}$ of $C$-heirs such that $f$ can be constructed from $E \cup \mathcal{E}$.
	\item[(2)]
	Let $f:C \to \mathbb{R}$ be purely $(e,r)$-prepared in $x$ with respect to $E$. Then $f$ can be constructed from $E$.
\end{itemize}

\vs{0.1cm}
{\bf Proof}

\vs{0.1cm}
Compare with Remark 3.11 in \cite{10} for a proof of (1). To obtain property (2) do the proof of Remark 3.11 in \cite{10} with $\mathcal{E}=\emptyset$. \hfill$\blacksquare$

\vs{0.4cm}
{\bf1.17 Fact} (Opris, \cite{10}, Theorem C)

\vs{0.1cm}
{\it
Let $e \in \mathbb{N}_0$. Let $X \subset \mathbb{R}^n \times \mathbb{R}$ be definable and let $E$ be a set of positive definable functions on $X$. Let $f:X \to \mathbb{R}$ be a function with exponential number at most $e$ with respect to $E$. Then there is $r \in \mathbb{N}_0$ and a definable cell decomposition $\mathcal{C}$ of $X_{\neq 0}$ such that for every $C \in \mathcal{C}$ there is a finite set $P$ of positive definable functions on $C$ and $\Theta:=(\Theta_0,...,\Theta_r)$ such that the function $f|_C$ is $(e,r)$-prepared in $x$ with center $\Theta$ with respect to $P$. Additionally the following holds.
\begin{itemize}
	\item [(1)]
	For every $g \in \log(P)$ there is $l \in \{-1,...,e-1\}$ such that $g$ is $(l,r)$-prepared in $x$ with center $\Theta$ with respect to $P$.
	\item [(2)]
	The following condition $(+_e)$ is satisfied: If $g \in \log(P)$ is $(l,r)$-prepared in $x$ with center $\Theta$ with respect to $P$ for $l \in \{-1,...,e-1\}$ then $g$ is a finite $\mathbb{Q}$-linear combination of functions from $\log(E)$ restricted to $C$ which have exponential number at most $l$ with respect to $E$.
\end{itemize}}

\section{Restricted Log-Exp-Analytic Functions}

\subsection{Basic Facts and Definitions}

The main results of this paper are formulated in the parametric setting. So we set up the concept of restricted log-exp-analytic functions in the variables $(t_1,...,t_n,x_1,...,x_m)$, where $x=(x_1,...,x_m)$ is serving as the tuple of independent variables of families of functions parameterized by $t:=(t_1,...,t_n)$. 

\vs{0.4cm}
Let $t$ range over $\mathbb{R}^n$ and $x$ over $\mathbb{R}^m$. We fix definable sets $C,X \subset \mathbb{R}^n \times \mathbb{R}^m$ with $C \subset X$. Suppose that $X_t$ is open for every $t \in \mathbb{R}^n$.

\vs{0.4cm}
{\bf2.1 Definition}

\vs{0.1cm}
Let $f:C \to \mathbb{R}$ be a function.
\begin{itemize}
	\item [(a)] We call $f$ \textbf{locally bounded in $x$ with reference set $X$} if the following holds. For $t \in \mathbb{R}^n$ and $w \in X_t$ there is an open neighbourhood $U$ of $w$ in $X_t$ such that $U \cap C_t = \emptyset$ or $f_t|_{U \cap C_t}$ is bounded. 
	\item [(b)] Suppose that $C_t$ is open for every $t \in \mathbb{R}^m$. We call $f$ \textbf{locally bounded in $x$} if $f$ is locally bounded in $x$ with reference set $C$.	
\end{itemize}

\vs{0.3cm}
{\bf2.2 Remark}

\vs{0.1cm}
Let $Y \subset \mathbb{R}^n \times \mathbb{R}^m$ be definable with $X \subset Y$ such that $Y_t$ is open for every $t \in \mathbb{R}^n$. Let $f:C \to \mathbb{R}$ be locally bounded in $x$ with reference set $Y$. Then $f$ is locally bounded in $x$ with reference set $X$.

\vs{0.5cm}
{\bf2.3 Remark}

\vs{0.1cm}
The set of locally bounded functions in $x$ with reference set $X$ on $C$ is an $\mathbb{R}$-algebra with respect to pointwise addition and multiplication. 

\vs{0.5cm}
{\bf2.4 Remark}

\vs{0.1cm}
Let $C_1,C_2 \subset X$ be disjoint definable sets such that $C = C_1 \cup C_2$. Let $g_j:C_j \to \mathbb{R}$ be locally bounded in $x$ with reference set $X$ for $j \in \{1,2\}$. Then 
$$g:C \to \mathbb{R}, (t,x) \mapsto \left\{\begin{array}{ll} g_1(t,x) , & (t,x) \in C_1,  \\
g_2(t,x), & (t,x) \in C_2, \end{array}\right.$$
is locally bounded in $x$ with reference set $X$.

\vs{0.5cm}
{\bf Proof}

\vs{0.1cm}
Let $t \in \mathbb{R}^n$ and let $w \in X_t$. Then there is an open neighbourhood $U_1$ of $w$ in $X_t$ such that $U_1 \cap C_1=\emptyset$ or $(g_1)_t|_{U_1 \cap C_1}$ is bounded and an open neighbourhood $U_2$ of $w$ in $X_t$ such that $U_2 \cap C_2=\emptyset$ or $(g_2)_t|_{U_2 \cap C_2}$ is bounded. Let $U:=U_1 \cap U_2$. We have that $U \cap C=\emptyset$ or by the definition of $g$ that $g_t|_{U \cap C}$ is bounded. 
\hfill$\blacksquare$

\vs{0.5cm}
{\bf2.5 Definition}

\vs{0.1cm}
Let $f:C \to \mathbb{R}$ be a function.
\begin{itemize}
	\item [(a)]
	Let $e \in \mathbb{N}_0$. We say that $f$ is \textbf{restricted log-exp-analytic in $x$ of order (at most) $e$ with reference set $X$} if $f$ has exponential number (at most) $e$ with respect to a set $E$ of positive definable functions on $C$ such that every $g \in \log(E)$ is locally bounded in $x$ with reference set $X$.
	
	\item [(b)]
	We say that $f$ is \textbf{restricted log-exp-analytic in $x$ with reference set $X$} if $f$ can be constructed from a set $E$ of positive definable functions on $C$ such that every $g \in \log(E)$ is locally bounded in $x$ with reference set $X$.
\end{itemize}

\vs{0.5cm}
{\bf2.6 Remark}

\vs{0.1cm}
Let $e \in \mathbb{N}_0$. Let $Y \subset \mathbb{R}^n \times \mathbb{R}^m$ be definable with $X \subset Y$ such that $Y_t$ is open for every $t \in \mathbb{R}^n$. Let $f:C \to \mathbb{R}$ be restricted log-exp-analytic in $x$ of order at most $e$ with reference set $Y$. Then $f$ is restricted log-exp-analytic in $x$ of order at most $e$ with reference set $X$.

\vs{0.5cm}
{\bf2.7 Remark}

\vs{0.1cm}
Let $l \in \mathbb{N}$. For $j \in \{1,...,l\}$ let $f_j:C \to \mathbb{R}$ be a function which is restricted log-exp-analytic in $x$ with reference set $X$. Let $F:\mathbb{R}^l \to \mathbb{R}$ be log-analytic. Then 
$$C \to \mathbb{R}, x \mapsto F(f_1(x),...,f_l(x)),$$
is restricted log-exp-analytic in $x$ with reference set $X$.

\vs{0.3cm}
{\bf Proof}

\vs{0.1cm}
This follows from Remark 1.8(2) applied to a set $E$ of positive definable functions on $C$ such that $f$ can be constructed from $E$ and every $g \in \log(E)$ is locally bounded in $x$ with reference set $X$. \hfill$\blacksquare$

\vs{0.5cm}
{\bf2.8 Remark}

\vs{0.1cm}
Let $C_1,C_2 \subset \mathbb{R}^n \times \mathbb{R}^m$ be disjoint and definable with $C_1 \cup C_2 = C$. For $j \in \{1,2\}$ let $f_j:C_j \to \mathbb{R}$ be restricted log-exp-analytic in $x$ with reference set $X$. Then
$$f:C \to \mathbb{R}, (t,x) \mapsto \left\{\begin{array}{ll} f_1(t,x) , & (t,x) \in C_1,  \\
f_2(t,x), & (t,x) \in C_2, \end{array}\right.$$
is restricted log-exp-analytic in $x$ with reference set $X$.

\vs{0.4cm}
{\bf Proof}

\vs{0.1cm}
For $j \in \{1,2\}$ let $E_j$ be a set of positive definable functions on $C_j$ such that every $g_j \in \log(E_j)$ is locally bounded in $x$ with reference set $X$ and $f_j$ can be constructed from $E_j$. Let
$$E:=\{g:C \to \mathbb{R} \mid \textnormal{ $g$ is a function with } g|_{C_j} \in E_j \textnormal{ for } j \in \{1,2\}\}.$$
We see with Remark 2.4 that every $g \in \log(E)$ is locally bounded in $x$ with reference set $X$. By Remark 1.7 $f$ can be constructed from $E$. \hfill$\blacksquare$

\vs{0.5cm}
{\bf2.9 Definition}

\vs{0.1cm}
A function $f:X \to \mathbb{R}$ is called \textbf{restricted log-exp-analytic in $x$} if $f$ is restricted log-exp-analytic in $x$ with reference set $X$.

\vs{0.5cm}
{\bf2.10 Remark}

\vs{0.1cm}
Let $k \in \mathbb{N}_0$. Let $w:=(w_1,...,w_k)$ range over $\mathbb{R}^k$. Let $g:\mathbb{R}^k \to \mathbb{R}^m$ be log-analytic and continuous. Let
$$V:=\{(t,x,w) \in X \times \mathbb{R}^k \mid (t,x+g(w)) \in X\}.$$
Let $f:X \to \mathbb{R}, (t,x) \mapsto f(t,x),$ be restricted log-exp-analytic in $x$. 
Then $F:V \to \mathbb{R}, (t,x,w) \mapsto f(t,x+g(w)),$ is restricted log-exp-analytic in $(x,w)$. 

\vs{0.4cm}
{\bf Proof}

\vs{0.1cm}
Note that $V_t$ is open in $\mathbb{R}^m \times \mathbb{R}^k$ for every $t \in \mathbb{R}^n$. Let $E$ be a set of positive definable functions on $X$ such that every $h \in \log(E)$ is locally bounded in $x$ and $f$ can be constructed from $E$. Consider
$$\tilde{E}:=\{V \to \mathbb{R}_{>0}, (t,x,w) \mapsto h(t,x+g(w)) \mid h \in E\}.$$  
Note that $F$ can be constructed from $\tilde{E}$. \hfill$\blacksquare$ 

\vs{0.4cm}
{\bf Claim}

\vs{0.1cm}
Let $h:X \to \mathbb{R}$ be locally bounded in $x$. Then $h^*:V \to \mathbb{R},(t,x,w) \mapsto h(t,x+g(w)),$ is locally bounded in $(x,w)$.

\vs{0.4cm}
{\bf Proof of the claim}

\vs{0.1cm}
Let $t \in \mathbb{R}^n$ and let $(x_0,w_0) \in V_t$. Note that $x_0+g(w_0) \in X_t$. Then there is an open neighbourhood $\hat{U}$ of $x_0+g(w_0)$ in $X_t$ such that $h_t|_{\hat{U}}$ is bounded. Let  
$$U:=\{(x,w) \in X_t \times \mathbb{R}^k \mid x+g(w) \in \hat{U}\}.$$
Since $g$ is continuous we see that $U$ is an open neighbourhood of $(x_0,w_0)$ in $V_t$ such that $(h^*)_t|_U$ is bounded.
\hfill$\blacksquare_{\textnormal{Claim}}$

\vs{0.3cm}
By the claim we see that $\tilde{E}$ is a set of locally bounded functions in $(x,w)$ on $V$ and we are done. \hfill$\blacksquare$

\subsection{A Preparation Theorem for Restricted \\
	Log-Exp-Analytic Functions}

\vs{0.3cm}
Let $m,l \in \mathbb{N}_0$ be with $n=l+m$. Let $w:=(w_1,...,w_l)$ range over $\mathbb{R}^l$ and $u:=(u_1,...,u_m)$ over $\mathbb{R}^m$. Let $C,X \subset \mathbb{R}^l \times \mathbb{R}^m \times \mathbb{R}$ be definable sets with $C \subset X$. Assume that $X_w$ is open for every $w \in \mathbb{R}^l$. 

\vs{0.5cm}
{\bf2.11 Definition}

\vs{0.1cm}
Let $e \in \mathbb{N}_0 \cup \{-1\}$ and $r \in \mathbb{N}_0$. We call $f:C \to \mathbb{R}, (w,u,x) \mapsto f(w,u,x),$ \textbf{$(m+1,X)$-restricted $(e,r)$-prepared in $x$} if $f$ is $(e,r)$-prepared in $x$ with respect to a set $E$ of positive definable functions such that every $g \in \log(E)$ is locally bounded in $(u,x)$ with reference set $X$.

\vs{0.5cm}
{\bf2.12 Proposition}

\vs{0.1cm}
{\it
Let $e \in \mathbb{N}_0 \cup \{-1\}$. Let $f:X \to \mathbb{R}$ be a restricted log-exp-analytic function in $(u,x)$ of order at most $e$. Then there is $r \in \mathbb{N}_0$ and a definable cell decomposition $\mathcal{C}$ of $X_{\neq 0}$ such that $f|_C$ is $(m+1,X)$-restricted $(e,r)$-prepared in $x$ for every $C \in \mathcal{C}$.
}

\vs{0.3cm}
{\bf Proof}

\vs{0.1cm}
Let $E$ be a set of positive definable functions such that every $g \in \log(E)$ is locally bounded in $(u,x)$ and $f$ has exponential number at most $e$ with respect to $E$. By Fact 1.17 there is a decomposition $\mathcal{C}$ of $X_{\neq 0}$ into finitely many definable cells such that the following holds for $C \in \mathcal{C}$. There is $r \in \mathbb{N}_0$ such that the function $f|_C$ is $(e,r)$-prepared with respect to a finite set $P$ of positive definable functions on $C$ such that every $g \in \log(P)$ is $(l,r)$-prepared in $x$ with respect to $P$ for $l \in \{0,...,e-1\}$ and if $g \in \log(P)$ is $(l,r)$-prepared in $x$ with respect to $P$ then $g$ is a finite $\mathbb{Q}$-linear combination of functions from $\log(E)$ with exponential number at most $l$ with respect to $E$ restricted to $C$ for $l \in \{0,...,e-1\}$. Thus by Remark 2.3 every $g \in \log(P)$ is locally bounded in $(u,x)$ with reference set $X$. \hfill$\blacksquare$

\vs{0.5cm}
{\bf2.13 Remark}

\vs{0.1cm}
Let $e \in \mathbb{N}_0 \cup \{-1\}$ and $r \in \mathbb{N}_0$. Let $f:C \to \mathbb{R}$ be $(m+1,X)$-restricted $(e,r)$-prepared in $x$. Then $f$ is not necessarily a restricted log-exp-analytic function in $(u,x)$ of order at most $e$ with reference set $X$. The assertion holds if $C=X$. 

\vs{0.4cm}
{\bf Proof}

\vs{0.1cm}
By Remark 1.16(1) there is a set $E$ of positive definable functions such that every $g \in \log(E)$ is locally bounded in $(u,x)$ with reference set $X$ and a set $\mathcal{E}$ of $C$-heirs such that $f$ can be constructed from $E \cup \mathcal{E}$. If $C=X$ one sees that $C_w$ is open for every $w \in \mathbb{R}^l$ and consequently that every $g \in \log(\mathcal{E})$ is locally bounded in $(u,x)$. (Let $r \in \mathbb{N}$ and $\Theta:=(\Theta_0,...,\Theta_r)$ be the center of a logarithmic scale on $C$. By an easy induction on $l \in \{0,...,r\}$ one can show with Definition 1.9(4) that $\Theta_l$ is locally bounded in $(u,x)$.) So $f$ is restricted log-exp-analytic in $(u,x)$. But this does not hold in general.

\vs{0.4cm}
Suppose that $w=0$ and $m=1$. Consider
$$C:=\{(u,x) \in \mathbb{R}^2 \mid 0 < u < 1, \tfrac{1}{1+u}+e^{-2/u+2e^{-1/u}} < x < \tfrac{1}{1+u}+e^{-1/u}\}.$$
 Let $X \subset \mathbb{R}^2$ be open with $C \subset X$ and $0 \in X$. Let
$$f:C \to \mathbb{R}, (u,x) \mapsto e^{-1/u}x.$$
An easy calculation shows that $C \to \mathbb{R}, u \mapsto e^{-1/u},$ is a $C$-heir (compare also with Example 2.37 in \cite{10}). So we see that $f$ is $0$-log-analytically prepared in $x$ with $C$-nice data. Let
$$
f^*:X \to \mathbb{R}, \textnormal{ } (u,x) \mapsto 
\left\{\begin{array}{ll} e^{-1/u}x,& (u,x) \in C, \\
0,& \textnormal{else}. \end{array}\right.
$$
Then for $x \in \textnormal{} ]0,1[$ we have that $f^* \sim e^{-1/u}x$ if $u \searrow 0$. So with Definition 2.5(a) one sees that $f$ is not restricted log-exp-analytic in $(u,x)$.
\hfill$\blacksquare$

\vs{0.5cm}
{\bf2.14 Definition}

\vs{0.1cm}
Let $e \in \mathbb{N}_0 \cup \{-1\}$ and $r \in \mathbb{N}_0$. We say that $f:C \to \mathbb{R}$ is \textbf{purely $(m+1,X)$-restricted $(e,r)$-prepared in $x$} if $f$ is purely $(e,r)$-prepared in $x$ with respect to a set $E$ of positive definable functions such that every $g \in \log(E)$ is locally bounded in $(u,x)$ with reference set $X$.

\vs{0.5cm}
{\bf2.15 Remark}

\vs{0.1cm}
Let $f:C \to \mathbb{R}$ be a function. Let $e \in \mathbb{N}_0 \cup \{-1\}$ and $r \in \mathbb{N}_0$. If $f$ is purely $(m+1,X)$-restricted $(e,r)$-prepared in $x$ then $f$ is restricted log-exp-analytic in $(u,x)$ with reference set $X$.

\vs{0.5cm}
{\bf2.16 Remark}

\vs{0.1cm}
Let $C \subset X$ be definable. Let $e \in \mathbb{N}_0 \cup \{-1\}$ and $r \in \mathbb{N}_0$. A function $f:C \to \mathbb{R}$ which is $(m+1,X)$-restricted $(e,r)$-prepared in $x$ is not purely $(m+1,X)$-restricted $(e,r)$-prepared in $x$ in general.

\vs{0.3cm}
{\bf Proof}

\vs{0.1cm}
This is easily seen with Remark 2.13 and Remark 2.15. \hfill$\blacksquare$

\section{Differentiability Properties of Restricted \\
	Log-Exp-Analytic Functions}

Let $\pi:\mathbb{R}^n \times \mathbb{R} \to \mathbb{R}^n, (t,x) \mapsto t,$ be the projection on the first $n$ coordinates.

\subsection{Simple Sets and Simple Preparation}

See also \cite{7} for Subsection 3.1.

\vs{0.4cm}
Let $C \subset \mathbb{R}^n \times \mathbb{R}$ be definable. Let $r \in \mathbb{N}_{0}$.

\vs{0.5cm} 
{\bf3.1 Definition}

\vs{0.1cm}
An $r$-logarithmic scale on $C$ is called \textbf{elementary} if its center is vanishing.

\vs{0.5cm}
{\bf3.2 Remark}

\vs{0.1cm}
An elementary $r$-logarithmic scale may not exist on $C$. If it exists it is uniquely determined and log-analytic.

\vs{0.5cm}
{\bf3.3 Definition}

\vs{0.1cm}
If $C$ has an elementary $r$-logarithmic scale then we call $C$ {\bf $r$-elementary}. The elementary $r$-logarithmic scale on $C$ is then denoted by $\ma{Y}_r^\mathrm{el}=\ma{Y}_{r,C}^\mathrm{el}$.

\vs{0.3cm}
For the next definition compare with the setting of \cite{5}, Section 4.

\vs{0.5cm}
{\bf3.4 Definition}

\vs{0.1cm}
We call $C$ {\bf simple} if for every $t \in \pi(C)$ we have $C_t=\,]0,d_t[$ for some $d_t \in \mathbb{R}_{>0} \cup \{+\infty\}$.

\vs{0.5cm}
{\bf3.5 Remark}

\vs{0.1cm}
Let $V:=\{(t,x) \in C \mid 0 \textnormal{ is interior point of }C_t\}$. Then $V$ is definable. Let $\mathcal{D}$ be a definable cell decomposition of $V_{\neq 0}$. Then
$$\pi(V)=\bigcup\{\pi(D) \mid D \in \mathcal{D} \textnormal{ simple}\}.$$

\vs{0.5cm}
We set $e_0:=0$ and $e_r:=\exp(e_{r-1})$ for $r\in \mathbb{N}$. In the following let $1/0:=\infty$.

\vs{0.5cm}
{\bf3.6 Fact} (Kaiser/Opris, \cite{7}, Section 2.2)

\vs{0.1cm}
{\it Let $C$ be simple and $r$-elementary and let $\ma{Y}_{r,C}^{\textnormal{el}}=(y_0,...,y_r)$. Then the following holds.
	\begin{itemize}
		\item[(1)] $\sup(C_t) \leq 1/e_r$ for all $t \in \pi(C)$. 
		\item[(2)] $y_0=x,y_1=\log(x), y_j=\log_{j-1}(-\log(x))$ for $j \in \{2,...,r\}$.
		\item[(3)] $\textnormal{sign}(\mathcal{Y}_{r,C}^{\textnormal{el}})=(1,-1,1,...,1) \in \mathbb{R}^{r+1}$.
\end{itemize}}

\vs{0.4cm}
{\bf3.7 Definition}

\vs{0.1cm}
We call $C$ \textbf{$r$-simple} if $C$ is simple and there is an $r$-logarithmic scale on $C$.

\vs{0.5cm}
{\bf3.8 Corollary}

\vs{0.1cm}
{\it Let $C$ be $r$-simple. Then an $r$-logarithmic scale on $C$ in uniquely determined.}

\vs{0.5cm}
{\bf3.9 Corollary}

\vs{0.1cm}
{\it Let $C$ be simple. Then the following are equivalent:
\begin{itemize}
	\item [(1)]
	$C$ is $r$-simple.
	\item [(2)] 
	$\sup C_t \leq 1/e_r$ for every $t \in \pi(C)$. 
\end{itemize}}

\vs{0.3cm}
With Fact 3.6 and Definition 3.7 one gets immediately the following (compare with the end of Section 2 in \cite{10}).

\vs{0.5cm}
{\bf3.10 Remark}

\vs{0.1cm}
{\it Let $C$ be simple. Then every $C$-nice function is log-analytic.}

\vs{0.3cm}
{\bf Proof}

\vs{0.1cm}
Let $g:\pi(C) \to \mathbb{R}$ be $C$-nice. Let $E$ be a set of $C$-heirs such that $g$ can be constructed from $E$. If $E=\emptyset$ we are clearly done. So suppose $E \neq \emptyset$. Let $h \in E$. There is $\hat{r} \in \mathbb{N}_0$, an $\hat{r}$-logarithmic scale $\hat{\mathcal{Y}}$ with center $(\hat{\Theta}_0,...,\hat{\Theta}_{\hat{r}})$ on $C$ and $l \in \{1,...,\hat{r}\}$ such that $h=\exp(\hat{\Theta}_l)$. Note that $\hat{\Theta}_l=0$ by Corollary 3.8. So we have $h=1$. We obtain $E=\{1\}$. By an easy indction on the exponential number $e \in \mathbb{N}_0$ of $g$ with respect to $E$ one sees that $g$ is log-analytic. \hfill$\blacksquare$

\vs{0.5cm}
{\bf3.11 Corollary}

\vs{0.1cm}
{\it Let $X \subset \mathbb{R}^n \times \mathbb{R}$ be definable and let $f:X \to \mathbb{R}$ be log-analytic of order at most $r$. Then there is a definable cell decomposition $\mathcal{C}$ of $X_{\neq 0}$ such that for every simple $C \in \mathcal{C}$ the cell $C$ is $r$-simple and $f|_C$ is $r$-log-analytically prepared in $x$ with center $0$, log-analytic coefficient and base functions.}

\vs{0.3cm}
{\bf Proof}

\vs{0.1cm}
By Fact 1.14 there is a definable cell decomposition $\mathcal{C}$ of $X_{\neq 0}$ such that for every $C \in \mathcal{C}$ we have that $f|_C$ is $r$-log-analytically prepared in $x$ with $C$-nice coefficient, $C$-nice base functions and $C$-nice center $\Theta$. Fix a simple $C \in \mathcal{C}$ and let 
$$(r,\mathcal{Y},a,q,s,v,b,P)$$
be a nice LA-preparing tuple for $f$. By Remark 3.10 we obtain that $a,b$ are log-analytic. Note also that $\Theta=0$. \hfill$\blacksquare$

\vs{0.3cm}
Compare also with \cite{7}, Section 2.3 for a proof of Corollary 3.11 explicitely for log-analytic functions on simple sets.

\vs{0.5cm}
{\bf3.12 Definition}

\vs{0.1cm}
Let $q=(q_0,...,q_r)\in \mathbb{Q}^{r+1}$ with $q\neq 0$. We set
$j(q):=\min\{j\mid q_j\neq 0\}$ and $\sigma(q):=\mathrm{sign}(q_{j(q)})\in \{\pm 1\}$. 
Moreover, let
$$q_\mathrm{diff}:=(q_0-1,\ldots,q_{j(q)}-1,q_{j(q)+1},\ldots,q_r).$$ 

\vs{0.5cm}
{\bf3.13 Remark}

\vs{0.1cm}
Let $C$ be $r$-simple. Let 
$q\in \IQ^{r+1}$ with $q\neq 0$.
Then
$$\lim_{x\searrow 0}|\ma{Y}_{r,C}^\mathrm{el}(x)|^{\otimes q}=
\left\{\begin{array}{cc}
0,&j(q)=0,\sigma(q)=+1,\\
+\infty,&j(q)=0,\sigma(q)=-1,\\
+\infty,&j(q)>0,\sigma(q)=+1,\\
0,&j(q)>0,\sigma(q)=-1.\\
\end{array}\right.$$

\vs{0.5cm}
{\bf3.14 Proposition}

\vs{0.1cm}
{\it Let $C$ be $r$-simple. Let $q\in \IQ^{r+1}$ with $q\neq 0$. Then
	$$\lim_{x \searrow 0} \Bigl|\frac{\frac{d}{dx}|\ma{Y}_{r,C}^\mathrm{el}(x)|^{\otimes q}}
	{|\ma{Y}_{r,C}^\mathrm{el}(x)|^{\otimes q_\mathrm{diff}}} \Bigl|=q_{j(q)}.$$}

\vs{0.1cm}
{\bf Proof}

\vs{0.1cm}
Compare with Kaiser/Opris in \cite{7}, Section 2.2. \hfill$\blacksquare$

\subsection{A Preparation Theorem for Restricted Log-Exp-Analytic Functions on Simple Sets}
Let $l,m \in \mathbb{N}_0$ be with $n=l+m$. Let $w$ range over $\mathbb{R}^l$ and $u$ over $\mathbb{R}^m$. Let $\pi_l:\mathbb{R}^l \times \mathbb{R}^m \to \mathbb{R}^l ,(w,u) \mapsto w$. Let $X \subset \mathbb{R}^n \times \mathbb{R}$ be definable such that $X_w$ is open in $\mathbb{R}^m \times \mathbb{R}$ for every $w \in \mathbb{R}^l$. Note that $\pi(X)_w$ is open in $\mathbb{R}^m$ for $w \in \mathbb{R}^l$. Assume that $0$ is an interior point of $X_t$ for every $t \in \pi(X)$.

\vs{0.3cm}
The goal for this section is to show that a restricted log-exp-analytic function $f:X \to \mathbb{R}, (w,u,x) \mapsto f(w,u,x),$ in $(u,x)$ can be log-analytically prepared in $x$ with restricted log-exp-analytic coefficient and base functions in $u$ with reference set $\pi(X)$ on simple sets. At first we prove that preparation implies pure preparation on simple cells.

\vs{0.5cm}
{\bf3.15 Remark}

\vs{0.1cm}
{\it 
Let $f:X \to \mathbb{R}$ be a restricted log-exp-analytic function in $(u,x)$. Then there are $r \in \mathbb{N}_0$, $e \in \mathbb{N}_0 \cup \{-1\}$ and a definable cell decomposition $\mathcal{C}$ of $X$ such that for every simple $C \in \mathcal{C}$ the cell $C$ is $r$-simple and $f|_C$ is purely $(m+1,X)$-restricted $(e,r)$-prepared in $x$.}

\vs{0.4cm}
{\bf Proof}

\vs{0.1cm}
By Proposition 2.12 there are $r \in \mathbb{N}_0$, $e \in \mathbb{N}_0 \cup \{-1\}$ and a definable cell decomposition $\mathcal{C}$ of $X$ such that for every $C \in \mathcal{C}$ the function $f|_C$ is $(m+1,X)$-restricted $(e,r)$-prepared in $x$. Fix a simple $C \in \mathcal{C}$ and a set $E$ of locally bounded functions in $(u,x)$ with reference set $X$ such that $f|_C$ is $(e,r)$-prepared in $x$ with respect to $E$. Note that $C$ is $r$-simple. We show by induction on $l \in \{-1,...,e\}$ that if $g \in \log(E) \cup \{f\}$ is $(l,r)$-prepared in $x$ with respect to $E$, then $g$ is purely $(l,r)$-prepared in $x$ with respect to $E$. For $l=-1$ there is nothing to show. 

\vs{0.2cm}
$l-1 \to l$: Let 
$$(r,\mathcal{Y},a,\exp(c),q,s,v,b,\exp(d),P)$$
be a preparing tuple for $g$. Note that $a$ and $b$ are $C$-nice and that $\mathcal{Y}=\mathcal{Y}_{r,C}^{\textnormal{el}}$. Therefore $\Theta=0$. Additionally we have that $a$ and $b$ are log-analytic by Remark 3.10. Therefore $g$ is purely $(l,r)$-prepared in $x$ with respect to $E$ by the inductive hypothesis and we are done.
\hfill$\blacksquare$

\vs{0.5cm}
{\bf3.16 Proposition}

\vs{0.1cm}
{\it 
Let $f:X \to \mathbb{R}, (w,u,x) \mapsto f(w,u,x),$ be restricted log-exp-analytic in $(u,x)$. Then there is $r \in \mathbb{N}_0$ and a definable cell decomposition $\mathcal{C}$ of $X$ such that for every simple $C \in \mathcal{C}$ the following holds. The cell $C$ is $r$-simple and $f|_C$ is $r$-log-analytically prepared in $x$ with coefficient and base functions which are restricted log-exp-analytic in $u$ with reference set $\pi(X)$.}

\vs{0.4cm}
{\bf Proof}

\vs{0.1cm}
By Remark 3.15 there are $r \in \mathbb{N}_0$, $e \in \mathbb{N}_0 \cup \{-1\}$ and a definable cell decomposition $\mathcal{Q}$ of $X$ such that for every simple $Q \in \mathcal{Q}$ the cell $Q$ is $r$-simple and $f|_Q$ is purely $(m+1,X)$-restricted $(e,r)$-prepared in $x$. Fix such a simple $Q \in \mathcal{Q}$. Let $d_y:=\sup(Q_y) \in \mathbb{R}_{>0} \cup \{\infty\}$ for $y \in \pi(Q)$. Let $E$ be a set of positive definable functions on $Q$ such that every $g \in \log(E)$ is locally bounded in $(u,x)$ with reference set $X$ and $f|_Q$ is purely $(e,r)$-prepared in $x$ with respect to $E$. We need the following claim.

\vs{0.4cm}
{\bf Claim}

\vs{0.1cm}
Let $h \in \log(E)$ be $r$-log-analytically prepared in $x$ with coefficient and base functions which are restricted log-exp-analytic in $u$ with reference set $\pi(X)$. Then there is a definable simple set $D \subset Q$ with $\pi(D)=\pi(Q)$ such that $h=h_1 + h_2$ where
\begin{itemize}
	\item [(1)] $h_1:\pi(D) \to \mathbb{R}$ is a function such that $\exp(h_1):\pi(D) \to \mathbb{R}$ is restricted log-exp-analytic in $u$ with reference set $\pi(X)$ and
	\item [(2)] $h_2:D \to \mathbb{R}$ is a bounded function such that $\exp(h_2)$ is $r$-log-analytically prepared in $x$ with coefficient $1$ and base functions which are restricted log-exp-analytic in $u$ with reference set $\pi(X)$.  
\end{itemize}

\vs{0.3cm}
{\bf Proof of the claim}

\vs{0.1cm}
Let
$$(r,\mathcal{Y},a,q,s,v,b,P)$$
be a corresponding LA-preparing tuple for $h$ where $b:=(b_1,...,b_s)$ and $P:=(p_1,...,p_s)^t$. Note that $\mathcal{Y}=\mathcal{Y}_{r,Q}^{\textnormal{el}}$. We have
$$h(t,x)=a(t) \vert{\mathcal{Y}(x)}\vert^{\otimes q} v(b_1(t)\vert{\mathcal{Y}(x)}\vert^{\otimes p_1},...,b_s(t)\vert{\mathcal{Y}(x)}\vert^{\otimes p_s})$$
for every $(t,x) \in Q$. We may assume that $a \neq 0$. Note that $q_0>0$ or if $q_0=0$ we have $q_{j(q)} \leq 0$. (Otherwise we have $\lim_{x \searrow 0} \vert{h(w,u,x)}\vert = \infty$ for every $(w,u) \in \pi(Q)$ by Remark 3.13, but for every $(w,u) \in \pi(Q)$ we have $(u,0) \in X_w$ and therefore there is an open neighbourhood $U$ of $(u,0)$ in $X_w$ such that $h_w|_{U \cap Q_w}$ is bounded.)

\vs{0.1cm}
We also have $p_{i0}>0$ or if $p_{i0}=0$ then $p_{i j(p_i)}\leq 0$ for every $i \in \{1,...,s\}$ by Remark 3.13.

\vs{0.3cm}
{\bf Case 1:} $q \neq 0$.

\vs{0.2cm}
For $t \in \pi(Q)$ let 
$$\hat{d}_t:=\sup\{x \in \textnormal{} ]0,d_t/2[ \textnormal{} \mid  \vert{a(t)}\vert \vert{\mathcal{Y}(c)}\vert^{\otimes q} \leq 1 \textnormal{ for every } c \in \textnormal{} ]0,x[\}$$
and $D:=\{(t,x) \in \pi(Q) \times \mathbb{R} \mid x \in \textnormal{} ]0,\hat{d}_t[\}$. Note that $D$ is definable and that $\pi(Q)=\pi(D)$. Set $h_1:=0$ and $h_2:=h|_D$. Note that $h_2$ is bounded. By the construction of $D$, using the exponential series and composition of power series we see that $\exp(h_2)$ is $r$-log-analytically prepared in $x$ as desired.

\vs{0.4cm}
{\bf Case 2:} $q = 0$.

\vs{0.3cm}
{\bf Subclaim 1}

\vs{0.1cm}
The coefficient $a$ is locally bounded in $u$ with reference set $\pi(X)$.

\vs{0.3cm}
{\bf Proof of Subclaim 1}

\vs{0.1cm}
Let $\hat{\pi}:\mathbb{R}^{m+1} \to \mathbb{R}^m, (u,x) \mapsto u$. Let $(w,u) \in \pi(X)$. Note that $(w,u,0) \in X$. Take an open ball $U$ around $(u,0)$ in $X_w$ such that either $U \cap Q_w = \emptyset$ or $h_w|_{U \cap Q_w}$ is bounded. So either $U \cap Q_w = \emptyset$ or $a_w|_{U \cap Q_w}$ is bounded. Let $B:=\pi(Q)$. Note that $\hat{\pi}(Q_w)=B_w$.

\vs{0.1cm}
If the former holds then $\hat{\pi}(U) \cap B_w = \hat{\pi}(U) \cap \hat{\pi}(Q_w) = \emptyset$. (Let $u \in \hat{\pi}(U) \cap \hat{\pi}(Q_w)$. Since $U$ is an open ball in $X_w$ there is $\epsilon>0$ such that $]0,\epsilon[ \subset U_u$. So $]0,\min\{\epsilon,d_{(w,u)}\}[ \subset U_u \cap (Q_w)_u$ and therefore $U \cap Q_w \neq \emptyset$. This shows also that $\hat{\pi}(U \cap Q_w)=\hat{\pi}(U) \cap \hat{\pi}(Q_w)=\hat{\pi}(U) \cap B_w$.)

\vs{0.1cm}
If the latter holds then $a_w|_{\hat{\pi}(U \cap Q_w)}$ is bounded since $a$ depends only on $(w,u)$. Since $\hat{\pi}(U \cap Q_w) = \hat{\pi}(U) \cap B_w$ we see that $a_w|_{\hat{\pi}(U) \cap B_w}$ is bounded. So take $\hat{\pi}(U) \subset \pi(X)_w$ which is an open neighbourhood of $u$ in $\pi(X)_w$. \hfill$\blacksquare_{\textnormal{Subclaim }1}$ 

\vs{0.5cm}
Let $\sum_{\alpha \in \mathbb{N}^s_0} c_{\alpha} X^{\alpha}$ be the power series expansion of $v$. Let 
$$\Gamma_1:=\{\alpha \in \mathbb{N}_0^s \mid {^t}P\alpha = 0\}$$
and
$\Gamma_2:=\mathbb{N}_0^s \setminus \Gamma_1$. Set $v_1:=\sum_{\alpha \in \Gamma_1} c_{\alpha}X^{\alpha}$ and $v_2:=\sum_{\alpha \in \Gamma_2} c_{\alpha}X^{\alpha}$. 
For $l \in \{1,2\}$ let 
$$g_l:Q \to \mathbb{R}, (t,x) \mapsto a(t) v_l(b_1(t)\vert{\mathcal{Y}(x)}\vert^{\otimes p_1},...,b_s(t)\vert{\mathcal{Y}(x)}\vert^{\otimes p_s}).$$
Let 
$$\mathcal{S}:=\{j \in \{1,...,s\} \mid p_j \neq 0\}.$$

\vs{0.3cm}
{\bf Subclaim 2} 

\vs{0.1cm}
There is a restricted log-exp-analytic function $\Psi:\pi(Q) \to \mathbb{R}, (w,u) \mapsto \Psi(w,u),$ in $u$ with reference set $\pi(X)$ which is locally bounded in $u$ with reference set $\pi(X)$ such that $g_1(w,u,x)=\Psi(w,u)$ for every $(w,u,x) \in Q$.

\vs{0.4cm}
{\bf Proof of Subclaim 2} 

\vs{0.1cm}
For $j \in \mathcal{S}$ we have that $\alpha_j=0$ for every $\alpha \in \Gamma_1$. Let $\{1,...,s\} \setminus \mathcal{S} = \{l_1,...,l_\lambda\}$ where $\lambda \in \{0,...,s\}$ and $l_1,...,l_\lambda \in \{1,...,s\}$ (if $\lambda=0$ then $\mathcal{S}=\{1,...,s\}$ and if $\lambda=s$ then $\mathcal{S}=\emptyset$). Note that there is a power series $\hat{v}$ which converges absolutely on an open neighbourhood of $[-1,1]^\lambda$ such that 
$$v_1(b_1(t)\vert{\mathcal{Y}(x)}\vert^{\otimes p_1},...,b_s(t)\vert{\mathcal{Y}(x)}\vert^{\otimes p_s})= \hat{v}(b_{l_1}(t),...,b_{l_\lambda}(t))$$
for $(t,x) \in Q$. Choose
$$\Psi:\pi(Q) \to \mathbb{R}, t \mapsto a(t) \hat{v}(b_{l_1}(t),...,b_{l_\lambda}(t)).$$
By Remark 2.7 $\Psi$ is restricted log-exp-analytic in $u$ with reference set $\pi(X)$ (since $\hat{v}$ defines a globally subanalytic function on $[-1,1]^\lambda$). Note that $\Psi$ has the desired properties since $a$ is locally bounded in $u$ with reference set $\pi(X)$. 

\hfill$\blacksquare_{\textnormal{Subclaim }2}$

\vs{0.4cm}
Let $\Psi$ be as in Subclaim 2. If $\mathcal{S}=\emptyset$ then $g_2=0$ and we are done by taking $D=Q$, $h_1:=\Psi$ and $h_2:=0$. So assume $\mathcal{S} \neq \emptyset$. Since $v$ is a power series which converges absolutely on an open neighbourhood of $[-1,1]^s$ there is $L \in \mathbb{R}_{>0}$ such that $\sum_{\alpha \in \mathbb{N}^s}\vert{c_\alpha}\vert < L$. Fix such an $L$. For $t \in \pi(Q)$ and $j \in \mathcal{S}$ let
$$\hat{d}_{j,t}:=\sup\{x \in \textnormal{}]0,d_t/2[ \textnormal{} \mid \vert{\mathcal{Y}(c)}\vert^{\otimes p_j} < \frac{1}{L\vert{a(t)b_j(t)}\vert} \textnormal{ for every }c \in \textnormal{} ]0,x[\}.$$
For $t \in \pi(Q)$ set $\hat{d}_t:=\min\{\hat{d}_{j,t} \mid j \in \mathcal{S}\}$. Consider
$$D:=\{(t,x) \in Q \mid x \in \textnormal{} ]0,\hat{d}_t[\}.$$
Note that $D$ is a simple definable set with $\pi(Q)=\pi(D)$.   

\vs{0.4cm}
{\bf Subclaim 3} 

\vs{0.1cm}
It is $\vert{g_2(t,x)}\vert \leq 1$ for every $(t,x) \in D$ and $\exp(g_2|_D)$ is $r$-log-analytically prepared in $x$ with coefficient $1$ and base functions which are restricted log-exp-analytic in $u$ with reference set $\pi(X)$.

\vs{0.4cm}
{\bf Proof of Subclaim 3}

\vs{0.1cm}
For $j \in \{1,...,s\}$ and $(t,x) \in D$ let $\phi_j(t,x):=b_j(t) \vert{\mathcal{Y}(x)}\vert^{\otimes p_j}$. For every $\alpha:=(\alpha_1,...,\alpha_s) \in \Gamma_2$ fix $i_\alpha \in \{1,...,s\}$ with $i_\alpha \in \mathcal{S}$ and $\alpha_{i_\alpha} \neq 0$. We have 
$$\vert{a(t)}\vert \vert{\phi_{i_\alpha}(t,x)}\vert \leq 1/L$$ 
for $(t,x) \in D$ and $\alpha \in \Gamma_2$. For $(t,x) \in D$ we obtain
\begin{eqnarray*}
	g_2(t,x) &=& a(t) v_2(\phi_1(t,x),...,\phi_s(t,x))\\
	&=& a(t)  \sum_{\alpha \in \Gamma_2}c_{\alpha} \phi_{i_{\alpha}}(t,x)^{\alpha_{i \alpha}} \prod_{j \neq i_{\alpha}} \phi_j(t,x)^{\alpha_j}\\
	&=&\sum_{\alpha \in \Gamma_2} c_{\alpha}a(t) \phi_{i_{\alpha}}(t,x) \phi_{i_{\alpha}}(t,x)^{\alpha_{i \alpha}-1} \prod_{j \neq i_{\alpha}} \phi_j(t,x)^{\alpha_j}\\
\end{eqnarray*}

and therefore $\vert{g_2(t,x)}\vert \leq 1/L \sum_{\alpha \in \Gamma_2} \vert{c_\alpha}\vert = 1$ since $\vert{\phi_j(t,x)}\vert \leq 1$ for $(t,x) \in D$ (because $D \subset Q$).

\vs{0.2cm}
Let $\mathcal{S}:=\{j_1,...,j_k\}$ where $k \in \mathbb{N}$. Note that $i_\alpha \in \{j_1,...,j_k\}$ is unique for $\alpha \in \mathbb{N}_0^s$. Let $(z_{j_1},...,z_{j_k})$ be a new tuple of variables. Consider for $l \in \{1,...,k\}$
$$\tilde{\phi}_l:D \to [-1,1], (t,x) \mapsto a(t)b_{j_l}(t)\vert{\mathcal{Y}(x)}\vert^{\otimes p_{j_l}}.$$
Consider 
$$\hat{v}:[-1,1]^{s+k} \to \mathbb{R}, (x_1,...,x_s,z_{j_1},...,z_{j_k}) \mapsto \sum_{\alpha \in \Gamma_2} c_{\alpha}z_{i_\alpha}x_{i_\alpha}^{\alpha_{i_\alpha}-1} \prod_{j \neq i_\alpha} x_j^{\alpha_j}.$$
Note that $\hat{v}$ is a well-defined globally subanalytic function since $\hat{v}$ defines a power-series which converges absolutely on an open neighbourhood of $[-1,1]^{s+k}$. Note that $\hat{v}([-1,1]^{s+k}) \subset [-L,L]$. We have for $(t,x) \in D$
$$g_2(t,x)=\hat{v}(\phi_1(t,x),...,\phi_s(t,x),\tilde{\phi}_1(t,x),...,\tilde{\phi}_k(t,x))$$
and therefore
$$\exp(g_2(t,x))=\exp^*(\hat{v}(\phi_1(t,x),...,\phi_s(t,x),\tilde{\phi}_1(t,x),...,\tilde{\phi}_k(t,x)))$$
where
$$\exp^*:\mathbb{R} \to \mathbb{R}, x \mapsto \left\{\begin{array}{ll} \exp(x), & x \in [-L,L], \\
	0, & \textnormal{else}. \end{array}\right.$$
is globally subanalytic. By using the exponential series and composition of power series we see that $\exp(g_2)$ has the desired properties (since $a$ and $b_j$ are restricted log-exp-analytic in $u$ with reference set $\pi(X)$). \hfill$\blacksquare_{\textnormal{Subclaim }3}$

\vs{0.4cm}
So take $h_1:=\Psi$ where $\Psi$ is as in Subclaim 2 and $h_2:=g_2|_D$. \hfill$\blacksquare_{\textnormal{Claim }}$

\vs{0.4cm}
Let $h \in \log(E) \cup \{f\}$ be purely $(l,r)$-prepared in $x$ with respect to $E$ where $l \in \{-1,...,e\}$. We show by induction on $l$ that there is a simple definable $A \subset Q$ with $\pi(A)=\pi(Q)$ such that $h|_A$ is $r$-log-analytically prepared in $x$ with coefficient and base functions which are restricted log-exp-analytic in $u$ with reference set $\pi(X)$. For $l=-1$ it is clear by choosing $A:=Q$. 

\vs{0.2cm}
$l-1 \to l:$ Let 
$$(r,\mathcal{Y},a,e^d,q,s,v,b,e^c,P)$$
be a purely preparing tuple for $h$ where $b:=(b_1,...,b_s)$, $e^c:=(e^{c_1},...,e^{c_s})$, and $P:=(p_1,...,p_s)^t$. Note that $\mathcal{Y}=\mathcal{Y}_{r,C}^{\textnormal{el}}$, that $a,b_1,...,b_s$ are log-analytic and that $d,c_1,...,c_s$ are purely $(l-1,e)$-prepared in $x$ with respect to $E$. We have
$$h(t,x)=a(t)\vert{\mathcal{Y}(x)}\vert^{\otimes q}e^{d(t,x)} v(b_1(t)\vert{\mathcal{Y}(x)}\vert^{\otimes p_1}e^{c_1(t,x)},...,b_s(t)\vert{\mathcal{Y}(x)}\vert^{\otimes p_s}e^{c_s(t,x)})$$
for every $(t,x) \in Q$. By the inductive hypothesis and the claim we find a simple definable set $A \subset Q$ with $\pi(A)=\pi(Q)$ and functions $d_1,c_{11},...,c_{1s}:\pi(A) \to \mathbb{R}$ and $d_2,c_{21},...,c_{2s}:A \to \mathbb{R}$ with the following properties: 

\begin{itemize}
	\item [(1)] The functions $\exp(d_1)$ and $\exp(c_{11}),...,\exp(c_{1s})$ are restricted log-exp-analytic in $u$ with reference set $\pi(X)$,
	\item [(2)] the functions $\exp(d_2),\exp(c_{21}),...,\exp(c_{2s})$ are $r$-log-analytically prepared in $x$ with coefficient $1$ and base functions which are restricted log-exp-analytic in $u$ with reference set $\pi(X)$,
	\item [(3)] we have $d|_A=d_1+d_2$ and $c_j|_A=c_{1j}+c_{2j}$ for $j \in \{1,...,s\}$.
\end{itemize}

Since $a$ and $b_1,...,b_s$ are log-analytic we see that the functions 
$$\hat{a}:\pi(A) \to \mathbb{R}, (w,u) \mapsto a(w,u) \exp(d_1(w,u)),$$
and
$$\hat{b}_j:\pi(A) \to \mathbb{R}, (w,u) \mapsto b_j(w,u) \exp(c_{1j}(w,u)),$$
for $j \in \{1,...,s\}$ are restricted log-exp-analytic in $u$ with reference set $\pi(X)$. For $(w,u,x) \in A$ we have 
$$h(w,u,x)=\hat{a}(w,u)\vert{\mathcal{Y}(x)}\vert^{\otimes q}e^{d_2(w,u,x)} v(\hat{\phi}_1(w,u,x),...,\hat{\phi}_s(w,u,x))$$
where $\hat{\phi}_j(w,u,x):=\hat{b}_j(w,u)\vert{\mathcal{Y}(x)}\vert^{\otimes p_j} e^{c_{2j}(w,u,x)}$ for $(w,u,x) \in A$ and $j \in \{1,...,s\}$. By composition of power series we obtain the desired $r$-log-analytical preparation for $h$ in $x$.  

\vs{0.2cm}
So we find a simple definable set $\hat{C} \subset Q$ with $\pi(\hat{C})=\pi(Q)$ such that $f|_{\hat{C}}$ is $r$-log-analytically prepared in $x$ with coefficient and base functions which are restricted log-exp-analytic in $u$ with reference set $\pi(X)$. With the cell decomposition theorem applied to every such $\hat{C}$ we are done (compare with \cite{3}, Chapter 3).\hfill$\blacksquare$

\vs{0.5cm}
An immediate consequence from this observation is the following.

\vs{0.5cm}
{\bf3.17 Proposition}

\vs{0.1cm}
{\it
Let $f:X \to \mathbb{R}, (w,u,x) \mapsto f(w,u,x)$, be restricted log-exp-analytic in $(u,x)$. Assume that $\lim_{x \searrow 0} f(t,x) \in \mathbb{R}$ for every $t \in \pi(X)$. Then $$h:\pi(X) \to \mathbb{R}, (w,u) \mapsto \lim_{x \searrow 0}f(w,u,x),$$
is restricted log-exp-analytic in $u$.}

\vs{0.4cm}
{\bf Proof}

\vs{0.1cm}
By Proposition 3.16 there is $r \in \mathbb{N}_0$ and a definable cell decomposition $\mathcal{C}$ of $X$ such that for every simple $C \in \mathcal{C}$ the cell $C$ is $r$-simple, and $f|_C$ is $r$-log-analytically prepared in $x$ with coefficient and base functions which are restricted log-exp-analytic in $u$ with reference set $\pi(X)$. Let $C \in \mathcal{C}$ be such a simple cell. Set $g:=f|_C$ and let 
$$(r,\mathcal{Y},a,q,s,v,b,P)$$
be a corresponding LA-preparing tuple for $g$. Note that $\mathcal{Y}=\mathcal{Y}_{r,C}^{\textnormal{el}}$. Then 
$$g(w,u,x)=a(w,u)\vert{\mathcal{Y}(x)}\vert^{\otimes q}v(b_1(w,u)\vert{\mathcal{Y}(x)}\vert^{\otimes p_1},...,b_s(w,u)\vert{\mathcal{Y}(x)}\vert^{\otimes p_s})$$ 
for $(w,u,x) \in C$. By the assumption, Remark 3.13 and Definition 1.13 we see that 
$$A:\pi(C) \to \mathbb{R}, (w,u) \mapsto \lim_{x \searrow 0} a(w,u)\vert{\mathcal{Y}(x)}\vert^{\otimes q},$$
and, for $j \in \{1,...,s\}$, that
$$B_j:\pi(C) \to [-1,1], (w,u) \mapsto \lim_{x \searrow 0} b_j(w,u)\vert{\mathcal{Y}(x)}\vert^{\otimes p_j},$$
are well-defined restricted log-exp-analytic functions in $u$ with reference set $\pi(X)$. We obtain for $(w,u) \in \pi(C)$
$$h(w,u)=A(w,u)v(B_1(w,u),...,B_s(w,u)).$$
Hence $h|_{\pi(C)}$ is restricted log-exp-analytic in $u$ with reference set $\pi(X)$ by Remark 2.7. So by Remark 2.8 we obtain that $h$ is restricted log-exp-analytic in $u$. \hfill$\blacksquare$

\vs{0.2cm}
Now we obtain Theorem A.

\vs{0.4cm}
{\bf Theorem A}

\vs{0.1cm}
{\it Let $Y \subset \mathbb{R}^l \times \mathbb{R}^m$ be definable such that $Y_w$ is open for every $w \in \mathbb{R}^l$. Let $f:Y \to \mathbb{R}, (w,u) \mapsto f(w,u),$ be restricted log-exp-analytic in $u$. Let $j \in \{1,...,m\}$ be such that $f_w$ is differentiable with respect to $u_j$ on $Y_w$ for every $w \in \pi_l(Y)$. Then $\partial f/\partial u_j$ is restricted log-exp-analytic in $u$.}

\vs{0.3cm}
{\bf Proof}

\vs{0.1cm}
We may assume that $f$ is differentiable with respect to the last variable $u_m$. We have to show that $\partial f/\partial u_m$ is restricted log-exp-analytic in $u$. Let $\textbf{e}_m:=(0,...,0,1) \in \mathbb{R}^m$ be the $m^{\textnormal{th}}$ unit vector. We define
$$V:=\{(w,u,x) \in X \times \mathbb{R} \mid (w,u+x\textbf{e}_m) \in X\}.$$
Note that $V_w$ is open and that $0$ is an interior point of $(V_w)_u$ for every $(w,u) \in \pi(V)$. Let
$$F:V \to \mathbb{R}, (w,u,x) \mapsto \left\{\begin{array}{cc}
 \frac{f(w,u+x\textbf{e}_m)-f(w,u)}{x}, & x \neq 0,\\
0, & \textnormal{ else}. \end{array}\right.$$
By Remark 2.10 we see that $G_1:V \to \mathbb{R}, (w,u,x) \mapsto f(w,u+x\textbf{e}_m),$ and $G_2:V \to \mathbb{R}, (w,u,x) \mapsto f(w,u)$, are restricted log-exp-analytic in $(u,x)$. So by Remark 2.7 we see that $F$ is restricted log-exp-analytic in $(u,x)$. Since
$$\frac{\partial f}{\partial u_m}(w,u)=\lim_{x \searrow 0}F(w,u,x)$$
for $(w,u) \in X$ we obtain by Proposition 3.17 that $\frac{\partial f}{\partial u_m}$ is a restricted log-exp-analytic function in $u$ with reference set $\pi(V)=Y$. \hfill$\blacksquare$

\subsection{Strong Quasianalyticity and Real Analyticity}

\vs{0.2cm}
{\bf The Univariate Case}

\vs{0.1cm}
We start our considerations with restricted log-exp-analytic functions in one variable and transfer these results into the multivariate case in the next paragraph.

\vs{0.2cm}
Let $X \subset \mathbb{R}^n \times \mathbb{R}$ be definable such that $X_t$ is open and $0$ is an interior point of $X_t$ for every $t \in \pi(X)$.
 
\vs{0.5cm}
{\bf3.18 Definition}

\vs{0.1cm}
Let $N \in \mathbb{N}$. Let $U \subset \mathbb{R}^n$ be open and let $g:U \to \mathbb{R}$ be a function. Let $a \in U$. 
The function $g$ is called {\bf $N$-flat} at $a$ if $g$ is $C^N$ at $a$ and all partial derivatives of $g$ of order at most $N$ vanish in $a$. The function $g$ is called {\bf flat} at $a$ if $g$ is $C^\infty$ at $a$ and all partial derivatives of $g$ vanish in $a$.

\vs{0.5cm}
The fact that a restricted log-exp-analytic function in $x$ can be log-analytically prepared in $x$ on simple definable cells implies the following.

\vs{0.5cm}
{\bf3.19 Proposition}

\vs{0.1cm}
{\it Let $f:X \to \mathbb{R}$ be restricted log-exp-analytic in $x$. Then there is $N \in \mathbb{N}$ such that the following holds for every $t \in \pi(X)$. If $f(t,-)$ is $N$-flat at $x=0$ then $f(t,-)$ vanishes identically on an open interval around $0 \in \mathbb{R}$.}

\vs{0.4cm}
{\bf Proof}

\vs{0.2cm}
The proof is the same as in \cite{7} for the log-analytic case. For the reader's convenience we give the details. 

\vs{0.1cm}
By also considering the function $f(t,-x)$ it is enough to show that the following holds:
There is $N \in \mathbb{N}$ such that for every $t \in \pi(X)$ with $f(t,-)$ is $N$-flat at $x=0$ we have $f(t,x)=0$ for all $x \in \, ]0,\varepsilon_t[$ for some $\epsilon_t \in \mathbb{R}_{>0}$. 
By Proposition 3.16 we find $r \in \mathbb{N}_0$ and a definable cell decomposition $\mathcal{C}$ of $X$ such that for every simple $C \in \mathcal{C}$ the following holds: $C$ is $r$-simple and $f|_C$ is $r$-log-analytically prepared in $x$. Fix a simple $C \in \ma{C}$ and let 
$$\big(r,\mathcal{Y},a,q,s,v,b,P \big)$$
be an LA-preparing tuple for $f|_C$ where $q:=(q_0,...,q_r)$. Note that $\mathcal{Y} = \mathcal{Y}_{r,C}^\textnormal{el}$. Choose $N_C \in \mathbb{N}$ such that $N_C > q_0$. Let $t \in \pi(C)$.
If $f(t,-)$ is $N_C$-flat at $x=0$ then $f(t,-)=o(x^{N_C})$ at $x=0$. But $|\mathcal{Y}(x)|^{\otimes q}/x^{N_C} \neq o(1)$ by Remark 3.13. Therefore we obtain $a(t)=0$. By Remark 3.5 we are done by taking 
$$N:=\max\{N_C \mid C \in \ma{C} \textnormal{ simple}\}.$$
\hfill$\blacksquare$

\vs{0.5cm}
{\bf3.20 Remark}

\vs{0.1cm}
Proposition 3.19 does not hold in $\mathbb{R}_{\textnormal{an,exp}}$ in general. Consider
$$f:\mathbb{R} \to \mathbb{R}, x \mapsto \left\{\begin{array}{cc}
e^{-\frac{1}{x}},&x>0,\\
0,& x \leq 0.
\end{array}\right.$$
Then $f$ is flat at $0$, but not the zero function. Note also that $f$ is $C^\infty$ at $0$, but not real analytic.

\vs{0.4cm}
Another consequence from Proposition 3.16 is the following statement about real analyticity of restricted log-exp-analytic functions in $x$ at $x=0$.

\vs{0.5cm}
{\bf3.21 Proposition}

\vs{0.1cm}
{\it Let $f:X \to \mathbb{R}$ be restricted log-exp-analytic in $x$. Then there is $N \in \mathbb{N}$ such that the following holds: If $f(t,-)$ is $C^N$ at $0$ then $f(t,-)$ is real analytic at $0$.}

\vs{0.4cm}
{\bf Proof}

\vs{0.1cm}
The proof is very similar to the proof in \cite{7} for the log-analytic case. By Proposition 3.16 we find $r \in \mathbb{N}_0$ and a definable cell decomposition $\mathcal{C}$ of $X$ such that for every simple $C \in \mathcal{C}$ the cell $C$ is $r$-simple and $f|_C$ is $r$-log-analytically prepared in $x$. Fix a simple cell $C \in \mathcal{C}$ and set $\eta_t:=\textnormal{sup }C_t$ for $t \in \pi(C)$. Let 
$$(r,\mathcal{Y},a,q,s,v,b,P)$$
be an LA-preparing tuple for $g:=f|_C$. Note that $\mathcal{Y}=\mathcal{Y}_{r,C}^{\textnormal{el}}$. Let $\sum_{\alpha\in \IN_0^s}c_\alpha X^\alpha$ be the power series expansion of $v$. Let 
$$\Gamma_1:=\big\{\alpha\in \IN_0^s\;\big\vert\;^tP\alpha+q\in \IN_0\times \{0\}^r\big\}$$
and
$\Gamma_2:=\IN_0^s\setminus \Gamma_1.$
Set
$v_1:=\sum_{\alpha\in \Gamma_1}c_\alpha X^\alpha$ and $v_2:=\sum_{\alpha\in \Gamma_2}c_\alpha X^\alpha.$
For $l\in \{1,2\}$ let
$$g_l:C\to \IR, (t,x)\mapsto a(t)|\ma{Y}(x)|^{\otimes q}v_l\big(b_1(t)|\ma{Y}(x)|^{\otimes p_1},\ldots,b_s(t)|\ma{Y}(x)|^{\otimes p_s}\big).$$
Then $g_1,g_2$ are restricted log-exp-analytic in $x$ with reference set $X$. It is $g=g_1+g_2$.
For $k\in \IN_0$ let
$$\Gamma_{1,k}:=\big\{\alpha\in \IN_0^s\;\big\vert\; ^tP\alpha+q=(k,0,\ldots,0)\big\}\subset \Gamma_1$$
and
$$d_k:\pi(C)\to \mathbb{R},t \mapsto
a(t)\sum_{\alpha\in \Gamma_{1,k}}c_\alpha\prod_{i=1}^sb_i(t)^{\alpha_i}.$$
Then
$$g_1(t,x)=\sum_{k=0}^\infty d_k(t)x^k$$
for $(t,x)\in C$.
The series to the right converges absolutely on $C$ and therefore $g_1$ extends to a well-defined extension
$$
\hat{g}_1:\hat{C}:=\big\{(t,x)\in \IR^{n+1}\;\big\vert\; t\in \pi(C), -\eta_t<x<\eta_t \big\}\to \IR,$$
$$(t,x)\mapsto \sum_{k=0}^\infty d_k(t)x^k,$$
such that $\hat{g}_1(t,-)$ is real analytic at $0$ for every $t \in \pi(C)$. Note that $\hat{C}$ is definable. By shrinking $\eta_t$ for $t \in \pi(X)$ if necessary we may assume that $\hat{C} \subset X$.

\vs{0.4cm}
{\bf Claim 1}

\vs{0.1cm}
The function $\hat{g}_1$ is restricted log-exp-analytic in $x$ with reference set $X$.

\vs{0.3cm}
{\bf Proof of Claim 1} 

\vs{0.1cm}
Note that $\hat{g}_1|_{\hat{C} \cap (\mathbb{R}^n \times \{0\})}=d_0$ is restricted log-exp-analytic in $x$ with reference set $X$ (since $d_0$ depends only on $t$). Let $D:=\hat{C} \cap (\mathbb{R}^n \times \mathbb{R}_{<0})$. We show that $\hat{g}_1|_D$ is restricted log-exp-analytic in $x$ with reference set $X$ and are done with Remark 2.8. Let $\Gamma_{1,e}:= \bigcup_{k \textnormal{ even}} \Gamma_{1,k}$ and $\Gamma_{1,o}:= \bigcup_{k \textnormal{ odd}} \Gamma_{1,k}$. Set $v_{1,e}:=\sum_{\alpha \in \Gamma_{1,e}}c_{\alpha}X^{\alpha}$ and $v_{1,o}:=\sum_{\alpha \in \Gamma_{1,o}}c_{\alpha}X^{\alpha}$. Then for $(t,x) \in \hat{C}$ with $x<0$ we have
$$\hat{g}_1(t,x)=a(t)\vert{\mathcal{Y}(-x)}\vert^{\otimes q} \Bigl(v_{1,e}(b_1(t)\vert{\mathcal{Y}(-x)}\vert^{\otimes p_1},...,b_s(t)\vert{\mathcal{Y}(-x)}\vert^{\otimes p_s})-$$
$$v_{1,o}(b_1(t)\vert{\mathcal{Y}(-x)}\vert^{\otimes p_1},...,b_s(t)\vert{\mathcal{Y}(-x)}\vert^{\otimes p_s}) \Bigl)$$  
which implies the desired assertion since it is easily seen that $\hat{g}_1$ can be constructed from a set $E$ of positive definable functions which depends only on $t$. (Note that $v_{1,e}$ and $v_{1,o}$ define globally subanalytic functions on $[-1,1]^s$. So choose $E$ such that $a$ and $b_1,...,b_s$ can be constructed from $E$.) \hfill$\blacksquare_{\textnormal{Claim }1}$

\vs{0.4cm}
Set
$$\hat{g}_2:\hat{C} \to \mathbb{R}, (t,x) \mapsto f(t,x) - \hat{g}_1(t,x).$$
Then $\hat{g}_2$ is restricted log-exp-analytic in $x$ with reference set $X$ by Claim 1. Note that $\hat{g}_2|_C=g_2$. Let $$\Lambda:=\big\{^tP\alpha+q\;\big\vert\;\alpha\in \Gamma_2 \big\}.$$
Then $\Lambda \subset \IQ^{r+1}\setminus (\IN_0\times \{0\}^r)$. Fix $t^* \in \pi(C)$. For $\lambda \in \Lambda$, let
$$\Gamma_{2,\lambda}:=\big\{\alpha\in \IN_0^s\;\big\vert\;^tP\alpha+q=\lambda\big\}$$
and 
$$e_{t^*,\lambda}:=a(t^*)\sum_{\alpha\in \Gamma_{2,\lambda}}c_\alpha\prod_{i=0}^sb_i(t^*)^{\alpha_i}.$$
Then 
$$\hat{g}_2(t^*,-)=\sum_{\lambda\in \Lambda}e_{t^*,\lambda}|\ma{Y}|^{\otimes \lambda}$$
on $]0,\eta_t[$. Let 
$$\Lambda_{t^*}:=\{\lambda \in \Lambda \mid e_{t^*,\lambda} \neq 0\}.$$
If $\Lambda_{t^*} = \emptyset$ then $\hat{g}_2(t^*,-)=0$ on $]0,\eta_{t^*}[$. If $\Lambda_{t^*} \neq \emptyset$ there is $\mu_{t^*} = (\mu_{t^*,0},...,\mu_{t^*,r}) \in \Lambda_{t^*}$ such that
$|\ma{Y}|^{\otimes \lambda}=o(|\ma{Y}|^{\otimes \mu_{t^*}})$ for all $\lambda\in \Lambda_{t^*}$ with $\lambda \neq \mu_{t^*}$.

\vs{0.4cm}
{\bf Claim 2}

\vs{0.1cm}
Let $M\in \IN$ be such that $f(t^*,-)$ is $C^M$ at $0$.
Then  $\mu_{t^*,0} \geq M$.

\vs{0.2cm}
{\bf Proof of Claim 2}

\vs{0.1cm}
Assume that $\mu_{t^*,0}<M$. 

\vs{0.2cm}
{\bf Case 1:} $\mu_{t^*,0}\in \mathbb{N}_0$. Then $m:=\mu_{t^*,0}+1\leq M$.
Differentiating $g_2$ $m$-times with respect to
$x$ we see with Proposition 3.14 that there is $\beta=(-1,\beta_1,\ldots,\beta_r)\in \mathbb{Q}^{r+1}$ such that
$$\lim_{x \searrow 0}\frac{\partial^{m}g_2/\partial x^{m}(t^*,x)}{|\ma{Y}|^{\otimes \beta}}\in \mathbb{R}^*.$$
Since $\hat{g}_2(t^*,-)=f(t^*,-) - \hat{g}_1(t^*,-)$ is $C^M$ at $0$ we obtain that
$$\lim\limits_{x \searrow 0} \frac{\partial^m\hat{g}_2}{\partial x^m}(t^*,x)=\frac{\partial^m\hat{g}_2}{\partial x^m}(t^*,0) \in \mathbb{R}$$
which contradicts that $\hat{g}_2(t^*,-)$ extends $g_2(t^*,-)$.

\vs{0.2cm}
{\bf Case 2:} $\mu_{t^*,0} \notin \mathbb{N}_0$. Then $m:=\lceil \mu_{t^*,0}\rceil \leq M$. Differentiating $g_2$ $m$-times with respect to $x$ we see with Proposition 3.14 (note that $\mu_{t^*} \neq 0$) that there is $\beta=(\beta_0,\beta_1,\ldots,\beta_r)\in \IQ^{r+1}$ with $\beta_0<0$ such that
$$\lim_{x \searrow 0}\frac{\partial^{m}g_2/\partial x^m(t^*,x)}{|\ma{Y}|^{\otimes \beta}}\in \IR^*.$$
But $\hat{g}_2(t^*,-)=f(t^*,-) - \hat{g}_1(t^*,-)$ is $C^M$ at $0$. We get the same contradiction as in Case $1$.
\hfill$\blacksquare_{\textnormal{Claim }2}$

\vs{0.4cm}
{\bf Claim 3}

\vs{0.1cm}
Let $M \in \mathbb{N}$ be such that $f(t^*,-)$ is $C^M$ at $0$. Then $\hat{g}_2(t^*,-)$ is $(M-1)$-flat at $0$.

\vs{0.2cm}
{\bf Proof of Claim 3}

\vs{0.1cm}
\textbf{Case 1:} $\Lambda_{t^*}=\emptyset$.
Then $\hat{g}_2(t^*,-)=0$ on $]0,\eta_{t^*}[$ and we are clearly done.

\vs{0.2cm}
\textbf{Case 2:} $\Lambda_{t^*} \neq \emptyset$. By Claim 1 we obtain that $\mu_{t^*,0} \geq M$. Hence we obtain by Proposition 3.14 and Remark 3.13 that
$$\lim_{x \searrow 0}\frac{\partial^m g_2}{\partial x^m}(t^*,x)=0$$
for all $m\in \{0,...,M-1\}$. Since $\hat{g}_2(t^*,-)=f(t^*,-)-\hat{g}_1(t^*,-)$ is $C^M$ at $0$ and since $\hat{g}_2(t^*,-)$ extends $g_2(t^*,-)$ we are done. \hfill$\blacksquare_{\textnormal{Claim }3}$

\vs{0.4cm}
Since the function $\hat{g}_2$ is restricted log-exp-analytic in $x$ with reference set $X$ we find by Proposition 3.19 some $K_C \in \mathbb{N}$ such that the following holds for every $t \in \pi(C)$: If $\hat{g}_2(t,-)$ is $K_C$-flat at $x=0$ then $\hat{g}_2(t,-)$ vanishes identically on some open interval around $0$. Set $N_C:=K_C+1$. Assume that $f(t,-)$ is $C^{N_C}$ at $0$. Then by Claim $3$ $\hat{g}_2(t,-)$ is $K_C$-flat and hence by the above that $f(t,-)=\hat{g}_1(t,-)$ on some open interval around $0$. Since $\hat{g}_1(t,-)$ is real analytic at $0$ we get that $f(t,-)$ is real analytic at $0$. By Remark 3.5 we are done by taking 
$$N:=\max\{N_C\mid C\in \ma{C}\mbox{ simple}\}.$$
\hfill$\blacksquare$

\vs{0.3cm}
{\bf3.22 Corollary}

\vs{0.1cm}
{\it Let $f:X \to \mathbb{R}$ be restricted log-exp-analytic in $x$. Assume that $f(t,-)$ is real analytic at $0$ for every $t \in \pi(X)$. Then there is a definable cell decomposition $\mathcal{B}$ of $\pi(X)$ such that $B \to \mathbb{R}, t \mapsto d^k/dx^k f(t,0),$ is real analytic for every $B \in \mathcal{B}$ and all $k \in \mathbb{N}_0$.}

\vs{0.4cm}
{\bf Proof}

\vs{0.1cm}
Using the notation of the previous proof we have $f(t,x)=\hat{g}_1(t,x)$ for all $(t,x)\in C$ where $C$ is a simple cell of the constructed definable cell decomposition $\mathcal{C}$ of $X$. Since functions definable in $\mathbb{R}_{\textnormal{an,exp}}$ are piecewise real analytic we find a definable cell decomposition $\mathcal{D}$ of $\pi(C) \subset \mathbb{R}^n$ such that the coefficient $a$ and the base functions $b_1,\ldots,b_s$ are real analytic on every $D \in \mathcal{D}$. Hence on each $D \in \mathcal{D}$ the coefficients $d_k$ of $\hat{g}_1$ are real analytic for every $k \in \mathbb{N}_0$.
Since $d^k/dx^k f(t,0)=k! d_k(t)$ we are done by Remark 3.5.
\hfill$\blacksquare$

\vs{1cm}
{\bf The Multivariate Case}

\vs{0.3cm}
Let $m \in \mathbb{N}$. Let $u:=(u_1,...,u_m)$ and $v:=(v_1,...,v_m)$ range over $\mathbb{R}^m$ and $x$ over $\mathbb{R}$. Let $u':=(u_1,...,u_{m-1})$. Let $\pi_n:\mathbb{R}^n \times \mathbb{R}^m \to \mathbb{R}^n ,(t,u) \mapsto t$. Let $X \subset \mathbb{R}^n \times \mathbb{R}^m$ be definable such that $X_t$ is open for every $t \in \mathbb{R}^n$.

\vs{0.5cm}
With Proposition 3.19 and familiar connectivity arguments we can easily prove theorem B, a strong quarianalyticity result of a restricted log-exp-analytic function.

\vs{0.5cm}
{\bf Theorem B}

\vs{0.1cm}
{\it Assume that $X_t$ is open and connected for every $t \in \mathbb{R}^n$. Let $f:X \to \mathbb{R}, (t,u) \mapsto f(t,u),$ be restricted log-exp-analytic in $u$. Then there is $N \in \mathbb{N}$ with the following property: Let $t \in \pi_n(X)$. If $f(t,-)$ is $C^N$ and if there is $a \in X_t$ such that $f(t,-)$ is $N$-flat at $a$ then $f(t,-)$ vanishes identically.} 	

\vs{0.4cm}
{\bf Proof}

\vs{0.1cm}
We set $U_\epsilon(a):=\{u \in \mathbb{R}^m \mid \vert{u-a}\vert<\epsilon\}$ where $a \in \mathbb{R}^m$ and $\epsilon \in \mathbb{R}_{>0}$. Let $y$ range over $\mathbb{R}$. Let 
$$\pi^*:\mathbb{R}^n \times \mathbb{R}^m \times \mathbb{R}^m \times \mathbb{R} \times \mathbb{R} \to \mathbb{R}^n \times \mathbb{R}^m \times \mathbb{R}^m \times \mathbb{R}, (t,u,v,y,x) \mapsto (t,u,v,y),$$
be the projection on the first $(n+2m+1)$ coordinates. Let 
$$\tilde{V}:=\{(t,u,v,y,x) \in X \times \mathbb{R}^m \times \mathbb{R} \times \mathbb{R} \mid (t,u+(x-y)v) \in X\},$$
and
$$G:\tilde{V} \to \mathbb{R}, (t,u,v,y,x) \mapsto f(t,u+(x-y)v).$$
By Remark 2.10 one sees that $G$ is restricted log-exp-analytic in $(u,v,y,x)$. Furthermore let
$$V:=\{(t,u,v,y,x) \in \tilde{V} \mid 0 \in \tilde{V}_{(t,u,v,y)}\}$$
and $F:=G|_V$. Then $F$ is restricted log-exp-analytic in $(u,v,y,x)$ with reference set $V$. Since $V_{(t,u,v,y)}$ is open for $(t,u,v,y) \in X \times \mathbb{R}^m \times \mathbb{R}$ we see that $F$ is also restricted log-exp-analytic in $x$. By Proposition 3.19 there is $N \in \mathbb{N}$ such that the following holds for every $(t,u,v,y) \in \pi^*(V)$: If $F(t,u,v,y,-)$ is $N$-flat at $x=0$ then $F(t,u,v,y,-)$ vanishes identically on an open interval around $0 \in \mathbb{R}$. Let $t \in \pi_n(X)$ be such that $f(t,-)$ is $C^N$ on $X_t$. Note that then $F(t,u,v,y,-)$ is $C^N$ on $V_{(t,u,v,y)}$ for $(t,u,v,y) \in \pi^*(V)$. Let $a \in X_t$ be such that $f(t,-)$ is $N$-flat at $a$. We show that this implies that $f(t,-)$ vanishes identically and are done. We start with the following 

\vs{0.3cm}
{\bf Claim}

\vs{0.1cm}
There is $r \in \mathbb{R}_{>0}$ such that $f(t,-)$ vanishes identically on $U_r(a)$.

\vs{0.3cm}
{\bf Proof of the claim}

\vs{0.1cm}
Let $r \in \mathbb{R}_{>0}$ be such that $U_{4r}(a)\subset X_t$. Then 
$$W:=U_r(a) \times U_1(a) \times ]-r,r[ \times ]-r,r[ \subset V_t.$$
So for $v \in U_1(a)$ we have that $F(t,a,v,0,-)$ is $N$-flat at $x=0$. Then by the above $F(t,a,v,0,-)$ vanishes on some open interval around $0$. Fix such a $v \in U_1(a)$ and let $A_v$ be the set of all $y \in \textnormal{} ]-r,r[$ such that $F(t,a,v,y,-)$ vanishes identically on some open interval around $0$. Then $A_v \neq \emptyset$ since $0 \in A_v$ by the above. Then $A_v$ is open. Let $y \in \overline{A_v} \textnormal{} \cap \textnormal{} ]-r,r[$. Then $F(t,a,v,y,-)$ is $N$-flat at $x=0$. Hence by the above $F(t,a,v,y,-)$ vanishes identically on some open interval around $0$. Therefore $A_v$ is closed in $]-r,r[$. Since intervals are connected we obtain that $A_v= \textnormal{}]-r,r[$ and hence that $F(t,a,v,-,-)$ vanishes identically on $]-r,r[ \times ]-r,r[$. Since $v \in U_1(a)$ is arbitrary we get that $f(t,-)$ vanishes identically on $U_r(a)$.
\hfill$\blacksquare_\textnormal{Claim}$

\vs{0.3cm}
Let $W_t$ be the set of all $z \in X_t$ such that $f$ vanishes identically on some open ball around $z$.
Then $W_t \neq \emptyset$ since $a \in W_t$ by the claim. Clearly $W_t$ is open.
Let $\hat{u} \in \overline{W_t} \cap X_t$. Then $f(t,-)$ is $N$-flat at $\hat{u}$. Again by the claim we get that $\hat{u} \in W_t$.
Therefore $W_t$ is closed in $X_t$. Since $X_t$ is connected we obtain that $W_t=X_t$ and hence that $f(t,-)$ vanishes identically on $X_t$.
\hfill$\blacksquare$

\vs{0.5cm}
With our results on restricted log-exp-analytic functions above we can establish a parametric version of Tamm's theorem for this class of functions. For this we adapt the reasoning of Van den Dries and Miller as in \cite{5}, Section 5.

\vs{0.5cm}
{\bf3.23 Definition}

\vs{0.1cm}
Let $U \subset \mathbb{R}^n$ be open and let $g:U \to \mathbb{R}$ be a function. Let $k \in \mathbb{N}$. Then $g$ is called {\bf $k$-times Gateaux-differentiable} or $G^k$ at $a$ if $x \mapsto g(a+xu)$ is $C^k$ at $x=0$ for every $u \in \mathbb{R}^n$ and $\big(u \mapsto d^kg(a+xu)/dx^k\big)(0)$ is given by a homogeneous polynomial in $u$ of degree $k$. The function $g$ is called $G^\infty$ at $a$ if $g$ is $G^k$ at $a$ for every $k\in \mathbb{N}$. 

\vs{0.3cm}
The following holds.

\vs{0.5cm}
{\bf3.24 Fact} (Van den Dries/Miller, \cite{5}, Section 2)

\vs{0.1cm}
{\it Let $U \subset \mathbb{R}^n$ be open and let $g:U \to \mathbb{R}$ be a function. Let $a \in U$. Then the following are equivalent:
	\begin{itemize}
		\item[(i)] The function $g$ is real analytic at $a$.
		\item[(ii)] The function $g$ is $G^\infty$ at $a$ and there is $\epsilon\in \mathbb{R}_{>0}$ such that for every $v \in \mathbb{R}^n$ with $\vert{v}\vert \leq 1$ the function $x \mapsto g(a+xv)$ is defined and real analytic on $]-\epsilon,\epsilon[$. 
\end{itemize}}

\vs{0.4cm}
{\bf3.25 Proposition}

\vs{0.1cm}
{\it Let $f:X \to \mathbb{R}, (t,u)\mapsto f(t,u)$, be a restricted log-exp-analytic function in $u$. Then there is $N \in \mathbb{N}$ such that the following holds for every $(t,u) \in X$: If $f(t,-)$ is $G^N$ at $u$ then $f(t,-)$ is $G^\infty$ at $u$.}

\vs{0.3cm}
{\bf Proof}

\vs{0.1cm}
Let $v:=(v_1,...,v_m)$. Let 
$$V:=\{(t,u,v,x) \in X \times \mathbb{R}^m \times \mathbb{R} \mid (t,u+xv) \in X\}.$$
Note that $V_t \subset \mathbb{R}^m \times \mathbb{R}^m \times \mathbb{R}$ is open and that $0$ is an interior point of $V_{(t,u,v)}$ for every $(t,u,v) \in X \times \mathbb{R}^m$. By Remark 2.10 we have that
$$g: V \to \mathbb{R}, (t,u,v,x) \mapsto f(t,u+xv),$$
is restricted log-exp-analytic in $(u,v,x)$. By Proposition 3.21  there is $K \in \mathbb{N}$ such that the following holds for every $(t,u,v) \in X \times \mathbb{R}^m$: If $x \mapsto f(t,u+xv)$ is $C^K$ at $0$ then $x \mapsto f(t,u+xv)$ is real analytic at $0$. For $k \in \mathbb{N}$ let $W_k$ be the set of all $(t,u) \in X$ such that $x \mapsto f(t,u+xv)$ is $C^k$ at $0$ for every $v\in \IR^m$.
We define
$$\Phi_k:X \times \IR^m \to \IR,
(t,u,v) \mapsto \left\{\begin{array}{ccc}
\frac{d^kf(t,u+xv)}{dx^k}(0),&&(t,u)\in W_k,\\
&\mbox{if}&\\
1,&&(t,u)\notin W_k.\\
\end{array}\right.$$
Then we have $W_k=W_K$ for all $k \geq K$. By Corollary 3.22 we find a definable cell decomposition $\ma{D}$ of $X \times \mathbb{R}^m$ such that $\Phi_k|_D$ is real analytic for every $D \in \mathcal{D}$ and all $k \in \mathbb{N}$. Let $\pi_X:\mathbb{R}^{n+m} \times \mathbb{R}^m \to \mathbb{R}^{n+m}, (t,u,v) \mapsto (t,u)$. For a definable cell decomposition $\mathcal{P}$ of $X \times \mathbb{R}^m$ we set 
$$\pi_X(\mathcal{P}):=\{\pi_X(P) \mid P \in \mathcal{P}\}$$ 
which is a definable cell decomposition of $X$.


\vs{0.3cm}
{\bf Claim 1}

\vs{0.1cm}
There is a definable cell decomposition $\mathcal{C}$ of $X \times \mathbb{R}^m$ compatible with $\mathcal{D}$ such that for every $B \in \pi_X(\mathcal{C})$ there is a non-empty open ball $U$ in $\mathbb{R}^m$ and $C^o \in \mathcal{C}$ with $\pi_X(C^o)=B$ such that for every $z \in B$ we have $U \subset (C^o)_z$.

\vs{0.4cm}
{\bf Proof of Claim 1}

\vs{0.1cm}
Let $z$ range over $\mathbb{R}^{n+m}$ and let $\hat{\pi}_X:\mathbb{R}^{n+m} \times \mathbb{R}^{m-1} \to \mathbb{R}^{n+m}, (z,u') \mapsto z,$ and $\pi^*: \mathbb{R}^{n+m} \times \mathbb{R}^m \to \mathbb{R}^{n+m}  \times \mathbb{R}^{m-1}, (z,u) \mapsto (z,u')$. We do an induction on $m$. 

\vs{0.2cm}
$m=1$: Let $Q \in \pi_X(\mathcal{D})$. It is enough to find a definable cell decomposition $\hat{\mathcal{C}}$ of $Q$ such that for every $\hat{C} \in \hat{\mathcal{C}}$ there is $A \in \mathcal{D}$ with $\hat{C} \subset \pi_X(A)$ and a non-empty open interval $I$ in $\mathbb{R}$ such that $I \subset A_z$ for every $z \in \hat{C}$. Of course this induces the desired cell decomposition $\mathcal{C}$ of $Q \times \mathbb{R}$. Let 
$$\mathcal{U}(\mathcal{D}_Q):=\{(z,u) \in Q \times \mathbb{R} \mid \textnormal{ there is $A \in \mathcal{D}$ such that $(z,u) \in A$}$$
$$\textnormal{and $A_y$ is open for every $y \in \mathbb{R}^{n+m}$}\}.$$

\vs{0.2cm}
By uniform finiteness there are $k \in \mathbb{N}_0$ and definable continuous functions $\sigma_1,...,\sigma_k:Q \to \mathbb{R}$ such that 
$$\mathcal{U}(\mathcal{D}_Q)=\{(z,u) \in Q \times \mathbb{R} \mid u \neq \sigma_j(z) \textnormal{ for }j \in \{1,...,k\}\}.$$
Let $\epsilon>0$. Pick $w_1,...,w_{k+1} \in \mathbb{R}$ such that $\vert{w_i-w_j}\vert>\epsilon$ for $i \neq j$. For $l \in \{1,...,k+1\}$ let $U_l:= \textnormal{}]w_l-\epsilon,w_l+\epsilon[$ and let 
$$Q_l:=\{z \in Q \mid \sigma_j(z) \notin U_l \textnormal{ for every }j \in \{1,...,k\}\}.$$
Note that $Q_l$ is definable, that $U_l \subset \mathcal{U}(\mathcal{D}_Q)_z$ for every $z \in Q_l$ and that $Q=\bigcup_{l=1}^{k+1}Q_l$. For $l \in \{1,...,k+1\}$ there is $D^o \in \mathcal{D}$ with $\pi_X(D^o) = Q$ such that $U_l \subset (D^o)_z$ for every $z \in Q_l$ since $U_l$ is connected. By choosing $\hat{\mathcal{C}}$ in this way that for every $\hat{C} \in \hat{\mathcal{C}}$ there is $l \in \{1,...,k+1\}$ such that $\hat{C} \subset Q_l$ we are done. 

\vs{0.3cm}
$m-1 \to m$: Let $\pi:\mathbb{R}^n \times \mathbb{R}^{m-1} \times \mathbb{R} \to \mathbb{R}^n \times \mathbb{R}^{m-1}, (t,u',u_m) \mapsto (t,u')$. Note that $\pi(\mathcal{D}):=\{\pi(D) \mid D \in \mathcal{D}\}$ is a definable cell decomposition of $X \times \mathbb{R}^{m-1}$. With the base case applied to every $\hat{D} \in \pi(\mathcal{D})$ we find a definable cell decomposition $\mathcal{T}$ of $X \times \mathbb{R}^{m-1}$ compatible with $\pi(\mathcal{D})$ such that for every $T \in \mathcal{T}$ there is a non-empty open interval $I_T$ in $\mathbb{R}$ and $A \in \mathcal{D}$ with $T \subset \pi(A)$ such that $I_T \subset A_{(z,u')}$ for every $(z,u') \in T$. By the inductive hypothesis there is a definable cell decomposition $\mathcal{K}$ of $X \times \mathbb{R}^{m-1}$ compatible with $\mathcal{T}$ such that for every $B_K \in \hat{\pi}_X(\mathcal{K}):=\{\hat{\pi}_X(K) \mid K \in \mathcal{K}\}$ there is $K^o \in \mathcal{K}$ with $\hat{\pi}_X(K^o)=B_K$ and a non-empty open ball $U_K$ in $\mathbb{R}^{m-1}$ such that $U_K \subset (K^o)_z$ for all $z \in B_K$. For $D \in \mathcal{D}$ and $K \in \mathcal{K}$ with $K \subset \pi(D)$ let 
$$C_{K,D}:=\{(z,u',u_m) \in D \mid (z,u') \in K\}.$$
Fix $K \in \mathcal{K}$ with a corresponding $K^o$ for $\hat{\pi}_X(K)$ as above. For $T^* \in \mathcal{T}$ with $K^o \subset T^*$ fix a corresponding $I_{T^*}$ and $A \in \mathcal{D}$ with $I_{T^*} \subset A_{(z,u')}$ for $(z,u') \in T^*$. Note that $C_{K^o,A}$ is a definable cell with $\hat{\pi}_X(K)=\pi_X(C_{K,D})=\pi_X(C_{K^o,A})$. Consider $U:=U_K \times I_{T^*}$. Let $z \in \pi_X(C_{K,D})$. By construction we have $U \subset (C_{K^o,A})_z$. So we see that
$$\mathcal{C}:=\{C_{K,D} \mid K \in \mathcal{K}, D \in \mathcal{D} \textnormal{ with } K \subset \pi(D)\}$$
is a definable cell decomposition of $X \times \mathbb{R}^m$ with the desired properties.
\hfill$\blacksquare_{\textnormal{Claim }1}$ 

\vs{0.4cm}
Let $\mathcal{C}$ be as in Claim 1. Let $\mathcal{B}:=\pi_X(\mathcal{C})$ and for $B \in \mathcal{B}$ let $\mathcal{D}_B:=\{C \in \mathcal{C} \mid \pi_X(C)=B\}$. Since $\mathcal{C}$ refines $\mathcal{D}$ we have that $\Phi_k|_C$ is real analytic for every $k \in \mathbb{N}$ and $C \in \mathcal{C}$.

\vs{0.5cm}
{\bf Claim 2}

\vs{0.1cm}
Let $k \in \mathbb{N}$. There is a definable function $w_k:X \times \mathbb{R}^m \to \mathbb{R}, (t,u,v) \mapsto w_k(t,u,v),$ such that $w_k|_C$ is real analytic for every $C \in \mathcal{C}$ and the following is equivalent for every $(t,u) \in X$.
	\begin{itemize}
		\item[(i)] The function $f(t,-)$ is $G^k$ at $u$.
		\item[(ii)] It is $w_k(t,u,v)=0$ for every $v \in \mathbb{R}^m$. 
\end{itemize}

\vs{0.4cm}
{\bf Proof of Claim 2}
 
\vs{0.1cm}
Let $\nu(k)$ be the dimension of the real vector space of homogeneous real polynomials of degree $k$ in the variables $V:=(V_1,\ldots,V_m)$
and let $M_1(V),\ldots,M_{\nu(k)}(V)$ be the homogeneous monomials of degree $k$ in $V$. For $p_1,...,p_{\nu(k)} \in \mathbb{R}^m$ let 
$$A(p_1,...,p_{\nu(k)}):=\left(\begin{array}{cccc}
M_1(p_1)&\cdot&\cdot&M_{\nu(k)}(p_1)\\
\cdot&& &\cdot\\
\cdot&& &\cdot\\
M_1(p_{\nu(k)})&\cdot&\cdot&M_{\nu(k)}(p_{\nu(k)})\\
\end{array}\right)\in M\big(\nu(k) \times \nu(k),\mathbb{R}).$$

Note that for $p_1,...,p_{\nu(k)} \in \mathbb{R}^m$ and all $s:=(s_1,...,s_{\nu(k)})$ the linear system of equations 
$$\sum_{j=1}^{\nu(k)} \kappa_jM_j(p_l) = s_l \textnormal{ }(*)$$
where $l \in \{1,...,\nu(k)\}$ has a unique solution $(\kappa_1,...,\kappa_{\nu(k)})$ if 
$$\textnormal{det}(A(p_1,...,p_{\nu(k)})) \neq 0.$$
We have that
$$T:=\{(p_1,...,p_{\nu(k)}) \in (\mathbb{R}^m)^{\nu_k} \mid \textnormal{det}(A(p_1,...,p_{\nu(k)})) = 0\}$$
is an algebraic set. So $(\mathbb{R}^m)^{\nu_\kappa} \setminus T$ is Zariski open and therefore dense in $(\mathbb{R}^m)^{\nu_\kappa}$. For $B \in \mathcal{B}$ consider the following: Fix $C^o \in \mathcal{D}_B$ and a non-empty open ball $U_B$ such that $U_B \subset (C^o)_z$ for every $z \in B$. 
Note that $(U_B)^{\nu(k)} \not \subseteq T$. Therefore there are points $p_{k,1},\ldots,p_{k,\nu(k)} \in U_B$ and linear functions $a_j:\mathbb{R}^{\nu(k)} \to \mathbb{R}$ for $j \in \{1,...,\nu(k)\}$ such that for all $s:=(s_1,...,s_{\nu(k)})$
$$P_{k,B}(s,V):=\sum_{j=1}^{\nu(k)}a_j(s)M_j(V)\in \IR[V]$$
is the unique homogeneous polynomial of degree $k$ with
$P_{k,B}(s,p_{k,i})=s_i$ for all $i \in \{1,\ldots,\nu(k)\}$.
We have $p_{k,j} \in (C^o)_z$ for every $z \in B$. Set
$$\hat{w}_{k,B}:B \times\IR^m \to \mathbb{R}, (t,u,v) \mapsto P_{k,B}\big(\Phi_k(t,u,p_{k,1}),\ldots,\Phi_k(t,u,p_{k,\nu(k)}),v\big)$$
and
$$w_{k,B}:=\hat{w}_{k,B}-\Phi_k|_{B \times \mathbb{R}^m}.$$
We have by the choice of $p_{k,1},...,p_{k,\nu(k)}$ that $B \to \mathbb{R}, z \mapsto \Phi_k(z,p_{k,j}),$ is real analytic for $j \in \{1,...,\nu_k\}$ and therefore that $w_{k,B}|_C$ is real analytic for $C \in \mathcal{D}_B$. 

\vs{0.3cm}
Letting $w_k:X \to \mathbb{R}$ be the function defined as $w_k(z)=w_{k,B}(z)$ if $z \in B$. Note that $w_k$ is well-defined, definable and that $w_k|_C$ is real analytic for every $C \in \mathcal{C}$. We show that $w_k$ fulfills the remaining requirements for $k \in \mathbb{N}$.

\vs{0.3cm}
$i) \Rightarrow ii)$: Let $(t,u) \in X$ be such that $f(t,-)$ is $G^k$ at $u$. Then $(t,u) \in W_k$ and 
$v \mapsto \Phi_k(t,u,v)$ is a homogeneous polynomial of degree $k$. Let $B \in \mathcal{B}$ be with corresponding points $p_{k,1},...,p_{k,\nu(k)} \in \mathbb{R}^m$ such that $(t,u) \in B$. By the definition of $\hat{w}_{k,B}$ we have
$\hat{w}_{k,B}(t,u,p_{k,j})=\Phi_k(t,u,p_{k,j})$ for all $j \in \{1,\ldots,\nu(k)\}$.
By the uniqueness of $P_{k,B}$ we obtain that
$$\hat{w}_{k,B}(t,u,v)=P_{k,B}\big(\Phi_k(t,u,p_{k,1}),\ldots,\Phi_k(t,u,p_{k,\nu(k)}),v\big)=\Phi_k(t,u,v)$$
and therefore $w_{k,B}(t,u,v)=0$ for all $v\in \IR^m$. 

\vs{0.2cm}
$ii) \Rightarrow i)$: Let $(t,u)\in X$ be such that $w_k(t,u,v)=0$ for all $v \in \IR^m$. Let $B \in \mathcal{B}$ be such that $(t,u) \in B$. Then
$\Phi_k(t,u,v)=\hat{w}_{k,B}(t,u,v)$ for all $v \in \IR^m$ and therefore
$v\mapsto \Phi_k(t,u,v)$ is a homogeneous polynomial of degree $k$. Since $k \geq 1$ it is not constant.
Hence we get that $(x,u)\in W_k$ and consequently $f$ is $G^k$ at $u$.
\hfill$\blacksquare_{\textnormal{Claim }2}$

\vs{0.4cm}
Let $w_k$ be as in the claim for $k \in \mathbb{N}$. By Tougeron (compare with \cite{12}) we find for every $C \in \ma{C}$ some $N_C\in \IN$ such that 
$$\bigcap_{k\in \IN}\{(t,u,v) \in C \mid w_k(t,u,v)=0\}=\bigcap_{k \leq N_C}\{(t,u,v) \in C \mid w_k(t,u,v)=0\}.$$
(Compare also with \cite{5}, 
Section 1). Let $N:=\max\{N_C\mid C \in \ma{C}\}$. Then 
$$\bigcap_{k\in \IN}\{(t,u,v) \in X \times \mathbb{R}^m \mid w_k(t,u,v)=0\} =$$
$$= \bigcap_{k \leq N}\{(t,u,v) \in X \times \mathbb{R}^m \mid w_k(t,u,v)=0\}.$$
Hence for every $(t,u)\in X$ we have that $f(t,-)$ is $G^\infty$ at $u$ if and only if $w_k(t,u,v)=0$ for all $k\in \IN$ and all $v \in \IR^m$ if and only if $w_k(t,u,v)=0$ for all $k \leq N$ and every $v\in \IR^m$
if and only if $f(t,-)$ is $G^N$ at $u$.
\hfill$\blacksquare$

\vs{0.5cm}
\textbf{Comparison with the proof of Van den Dries and Miller:}

\vs{0.1cm}
There is a little gap in the proof in \cite{5} (and therefore also in \cite{7}, Section 3.3). Although Van den Dries and Miller worked with a definable cell decomposition $\mathcal{C}$ of $X \times \mathbb{R}^m$ such that $\Phi_k|_C$ is real analytic for $k \in \mathbb{N}$ and $C \in \mathcal{C}$, the points $p_1,...,p_{\nu(k)} \in \mathbb{R}^m$ are choosen ''arbitrarily'' and even independently of the decomposition $\mathcal{C}$ such that the linear system of equations $(*)$ above has a unique solution $\kappa$. So it is assumed without justification that $w_k|_C$ is real analytic for $C \in \mathcal{C}$ where
$$w_k(t,u,v):=P_k(\Phi_k(t,u,p_1),...,\Phi_k(t,u,p_{\nu(k)}),v)-\Phi_k(t,u,v)$$
for $k \in \mathbb{N}$ and $(t,u,v) \in X \times \mathbb{R}^m$ (i.e. that $\Phi_k(t,u,p_j)|_C$ for $j \in \{1,...,\nu(k)\}$ and therefore that $P_k(\Phi_k(t,u,p_1),...,\Phi_k(t,u,p_{\nu(k)}),v)|_C$ is real analytic for $C \in \mathcal{C}$).

\newpage
Finally we obtain Theorem C.

\vs{0.5cm}
{\bf Theorem C} 

\vs{0.1cm}
{\it Let $f:X \to \mathbb{R}, (t,u) \mapsto f(t,u),$ be a restricted log-exp-analytic function in $u$. Then there is $N \in \mathbb{N}$ such that for all $t \in \mathbb{R}^n$ if $f(t,-)$ is $C^N$ at $u$ then $f(t,-)$ is real analytic at $u$.} 

\vs{0.3cm}
{\bf Proof}

\vs{0.1cm}
By Proposition 3.25 there is $N_1 \in \mathbb{N}$ such that for every $(t,u) \in X$ if $f(t,-)$ is $G^{N_1}$ at $u$ then $f(t,-)$ is $G^{\infty}$ at $u$. Let 
$$V:=\{(t,u,v,x) \in X \times \mathbb{R}^m \times \mathbb{R} \mid (t,u+xv) \in X\}$$
and
$$F:V \to \mathbb{R}, (t,u,v,x)\mapsto f(t,u+xv).$$ 
Since $X_t$ is open we see that $V_t$ is open for every $t \in \mathbb{R}^n$. By Remark 2.10 $F$ is restricted log-exp-analytic in $(u,v,x)$. By Proposition 3.21 there is $N_2 \in \mathbb{N}$ such that if $F(t,u,v,-)$ is $C^{N_2}$ at $0$ then $F(t,u,v,-)$ is real analytic at $0$. Taking $N:=\max\{N_1,N_2\}$ we are done with Fact 3.24.
\hfill$\blacksquare$

\vs{0.5cm}
{\bf3.26 Corollary}

\vs{0.1cm}
{\it Let $f:X \to \mathbb{R}, (t,u) \mapsto f(t,u),$ be restricted log-exp-analytic in $u$. Then the set of all $(t,u) \in X$ such that $f(t,-)$ is real analytic at $u$ is definable.}

\vs{1cm}
Andre Opris\\
University of Passau\\
Faculty of Computer Science and Mathematics\\
andre.opris@uni-passau.de\\
D-94030 Germany

\end{document}